\documentclass[leqno,a4paper,12pt]{article}

\usepackage{amsmath,amsthm}

\newtheorem{Thm}{\indent Theorem}[section]
\newtheorem{Prop}[Thm]{\indent Proposition}
\newtheorem{Cor}[Thm]{\indent Corollary}
\newtheorem{Lemma}[Thm]{\indent Lemma}

\theoremstyle{definition}
\newtheorem{Rem}[Thm]{\indent Remark}
\newtheorem{Con}[Thm]{\indent Convention}
\newtheorem{Def}[Thm]{\indent Definition}

\def\qed{{\hskip0pt\unskip\unskip\nobreak\hfil\penalty50
          \hskip1em\hbox{}\nobreak\hfil
          {\bf q.e.d.}%
          \parfillskip=0pt\finalhyphendemerits=0
          \par}\medskip}

\newenvironment{Proof}
               {{\it Proof.}\quad}
               {\qed}

\newenvironment{Proofof}[1]
               {{\it Proof of #1.}\quad}
               {\qed}

\newcommand{\Prime}{\kern3\fontdimen1\font$'$\kern-7\fontdimen1\font}

\long\def\forget#1{}

\long\def\beginSIDEREMARK#1\endSIDEREMARK
    {{\par\bigskip\advance\leftskip by 2cm
                  \advance\rightskip by -2cm\noindent
      {\bf Our own side remark:} #1
      \par\bigskip\noindent}}

\long\def\beginFORGET#1\endFORGET{#1}
\long\def\beginFORGET#1\endFORGET{}

\def\?{\ ???\ \immediate\write16{}%
\immediate\write16{Warning: There was still a question mark . . . }%
\immediate\write16{}}

\usepackage{amsmath}

\usepackage{exscale}
\usepackage{amssymb}

\newcommand{\BA}{{\mathbb{A}}}

\newcommand{\BD}{{\mathbb{D}}}

\newcommand{\BG}{{\mathbb{G}}}

\newcommand{\BZ}{{\mathbb{Z}}}

\newcommand{\Fa}{{\mathfrak{a}}}

\newcommand{\FF}{{\mathfrak{F}}}

\newcommand{\FN}{{\mathfrak{N}}}

\newcommand{\FW}{{\mathfrak{W}}}
\newcommand{\FX}{{\mathfrak{X}}}
\newcommand{\FY}{{\mathfrak{Y}}}

\newcommand{\CD}{{\cal D}}

\newcommand{\CL}{{\cal L}}

\newcommand{\CO}{{\cal O}}

\newcommand{\CZ}{{\cal Z}}

\newcommand{\Ext}{\mathop{\rm Ext}\nolimits}
\newcommand{\coker}{\mathop{\rm Coker}\nolimits}
\newcommand{\loccit}{[loc.$\;$cit.]}

\def\tei{\, | \,}

\def\id{{\rm id}}

\newbox\mybox
\def\arrover#1{\mathrel{
       \setbox\mybox=\hbox spread 1.4em{\hfil$\scriptstyle#1$\hfil}
       \vbox{\offinterlineskip\copy\mybox
             \hbox to\wd\mybox{\rightarrowfill}}}}
\def\larrover#1{\mathrel{
       \setbox\mybox=\hbox spread 1.4em{\hfil$\scriptstyle#1$\hfil}
       \vbox{\offinterlineskip\copy\mybox
             \hbox to\wd\mybox{\leftarrowfill}}}}

\def\ontoover#1{\mathrel{
       \setbox\mybox=\hbox spread 1.4em{\hfil$\scriptstyle#1$\hfil}
       \vbox{\offinterlineskip\copy\mybox
             \hbox to\wd\mybox{\rightarrowfill\hskip-2.8mm
                               $\rightarrow$}}}}
\def\leftontoover#1{\mathrel{
       \setbox\mybox=\hbox spread 1.4em{\hfil$\scriptstyle#1$\hfil}
       \vbox{\offinterlineskip\copy\mybox
             \hbox to\wd\mybox{$\leftarrow$\hskip-2.8mm
                               \leftarrowfill}}}}
\def\longto{\longrightarrow}
\def\into{\hookrightarrow}

\def\isoto{\arrover{\sim}}

\def\longinto{\lhook\joinrel\longrightarrow}

\usepackage[curve,matrix,arrow,cmtip]{xy}
\NoComputerModernTips

\def\myxymessage{\def\messagetext
   {Here an xy-pic diagram was omitted to speed up compilation . . . }
   \immediate\write16{\messagetext}
   \hbox{\bf \messagetext}}
\def\filxymatrix#1{\myxymessage}
\def\filxyarray#1{\myxymessage}

\newdir^{ (}{{}*!/-3pt/\dir^{(}}
\newdir_{ (}{{}*!/-3pt/\dir_{(}}
\newdir^{ )}{{}*!/+3pt/\dir^{)}}
\newdir_{ )}{{}*!/+3pt/\dir_{)}}

\def\rscript#1{\hbox to 0pt{$\scriptstyle#1$\hss}}

\newcommand{\cone}{\mathop{\rm Cone}\nolimits}
\newcommand{\DM}{\mathop{DM^{eff}_-(k)}\nolimits}
\newcommand{\DeffgM}{\mathop{DM^{eff}_{gm}(k)}\nolimits}
\newcommand{\DgM}{\mathop{DM_{gm}(k)}\nolimits}
\newcommand{\Hom}{\mathop{\rm{Hom}}\nolimits}

\newcommand{\iastXmjZ}{\mathop{i_{(\doX)_m}^* j_* \BZ}\nolimits}
\newcommand{\iastXnmjZ}{\mathop{i_m^* j_* \BZ}\nolimits}
\newcommand{\iXjZ}{\mathop{i_{\doX}^! j_! \, \BZ}\nolimits}
\newcommand{\iXmjZ}{\mathop{i_{(\doX)_m}^! j_! \, \BZ}\nolimits}
\newcommand{\iXnmjZ}{\mathop{i_m^! j_! \, \BZ}\nolimits}
\newcommand{\iastYjDZ}{\mathop{i_Y^* j_* \BD(\BZ)}\nolimits}
\newcommand{\iastYmjZ}{\mathop{i_{Y_m}^* j_* \BZ}\nolimits}
\newcommand{\iastYjZ}{\mathop{i_Y^* j_* \BZ}\nolimits}
\newcommand{\iastYjaZ}{\mathop{i_Y^* j_{\alpha *} \BZ}\nolimits}
\newcommand{\iYjZ}{\mathop{i_Y^! j_! \, \BZ}\nolimits}
\newcommand{\iYmjZ}{\mathop{i_{Y_m}^! j_! \, \BZ}\nolimits}
\newcommand{\Mgm}{\mathop{M_{gm}}\nolimits}
\newcommand{\Mcgm}{\mathop{M_{gm}^c}\nolimits}
\newcommand{\dMgm}{\mathop{\partial M_{gm}}\nolimits}
\newcommand{\oT}{\mathop{\overline{T}}\nolimits}

\newcommand{\oV}{\mathop{\overline{V}}\nolimits}
\newcommand{\oW}{\mathop{\overline{W}}\nolimits}
\newcommand{\oX}{\mathop{\overline{X}}\nolimits}
\newcommand{\doX}{\mathop{\partial \overline{X}}\nolimits}
\newcommand{\RC}{\mathop{{\bf R} C}\nolimits}
\newcommand{\Schi}{\mathop{Sch^\infty \! /k}\nolimits}
\newcommand{\ShN}{\mathop{Shv_{Nis}(SmCor(k))}\nolimits}
\newcommand{\SmC}{\mathop{SmCor(k)}\nolimits}
\newcommand{\Spec}{\mathop{{\rm Spec}}\nolimits}
\newcommand{\uC}{\mathop{\underline{C}}\nolimits}
\newcommand{\uh}{\mathop{\underline{h}}\nolimits}
\newcommand{\uHom}{\mathop{\underline{Hom}}\nolimits}
\newcommand{\zeff}{\mathop{z_{equi}^{eff}}\nolimits}
\newcommand{\zeq}{\mathop{z_{equi}}\nolimits}

\let\oldbullet\bullet
\def\bullet{{\mathchoice{\oldbullet}%
                        {\oldbullet}%
                        {\scriptscriptstyle\oldbullet}%
                        {\oldbullet}}}

\newcommand{\argdot}{{\;\bullet\;}}

\begin{document}

%

\hfuzz=3pt
\overfullrule=10pt                   


\setlength{\abovedisplayskip}{6.0pt plus 3.0pt}
\setlength{\belowdisplayskip}{6.0pt plus 3.0pt}
\setlength{\abovedisplayshortskip}{6.0pt plus 3.0pt}
\setlength{\belowdisplayshortskip}{6.0pt plus 3.0pt}

\setlength{\baselineskip}{13.0pt}
\setlength{\lineskip}{0.0pt}
\setlength{\lineskiplimit}{0.0pt}

%
%

\title{The boundary motive: definition and basic properties
\footnotemark
\footnotetext{To appear in Compositio Mathematica.}
}
\author{\footnotesize by\\ \\
\mbox{\hskip-2cm
\begin{minipage}{6cm} \begin{center} \begin{tabular}{c}
J\"org Wildeshaus \\[0.2cm]
\footnotesize Institut Galil\'ee\\[-3pt]
\footnotesize Universit\'e Paris 13\\[-3pt]
\footnotesize Avenue Jean-Baptiste Cl\'ement\\[-3pt]
\footnotesize F-93430 Villetaneuse\\[-3pt]
\footnotesize France\\
{\footnotesize \tt wildesh@math.univ-paris13.fr}
\end{tabular} \end{center} \end{minipage}
\hskip-2cm}
\\[2.5cm]
}
\date{September 3, 2005}
\maketitle
\quad \\[-1.7cm]
\begin{abstract}
\noindent
We introduce the notion of the boundary motive of a scheme $X$
over a perfect field. By definition, it measures the difference
between the motive $\Mgm (X)$ and the motive with compact support
$\Mcgm (X)$, as defined and studied in \cite{VSF}. We develop
three tools to compute the boundary motive in terms of 
the geometry of a compactification
of $X$: co-localization, invariance under abstract blow-up,
and analy\-ti\-cal invariance. 
We then prove auto-duality of the boundary motive
of a smooth scheme $X$.
As a formal consequence of this, and of co-localization, we
obtain a fourth computational tool, namely
localization for the boundary motive.  
In a sequel to this work \cite{W9}, these tools will be applied to
Shimura varieties. \\


\end{abstract}


\bigskip
\bigskip
\bigskip

\noindent {\footnotesize Math.\ Subj.\ Class.\ (2000) numbers:
14F42 (14C25, 19E15).
}

\eject
\tableofcontents

\bigskip
\vspace*{0.5cm}


%
%

\setcounter{section}{-1}
\section{Introduction}
\label{Intro}



In this paper, we define the boundary motive 
$\dMgm (X)$ of a scheme $X$ over a perfect field $k$. 
By its very construction, it is part of an exact triangle
\[
\dMgm (X) \longto \Mgm (X) \longto \Mcgm (X) \longto \dMgm (X) [1] \; ,
\]
where $\Mgm (X)$ and $\Mcgm (X)$ denote the motive of $X$ and its
motive with compact support, respectively. The exact triangle is
in the triangulated category $\DM$ of effective motivic complexes
which, as  $\Mgm (X)$ and $\Mcgm (X)$ was defined in \cite{VSF}.
We refer to Section~\ref{1} for a review of these constructions. 
We expect this exact triangle to be of a certain interest. First, it
induces long exact sequences for motivic homology and cohomology. More
generally, any exact functor 
to a triangulated category $D$ will 
induce long exact sequences of $\Ext$-groups in $D$ which are a priori
compatible with the sequence in motivic homology resp.\ cohomology. 
Second, the exact triangle itself
can be employed to construct explicit
extensions of objects in $\DM$, i.e., classes in motivic cohomology. \\

Note that most of the existing attempts to prove the Beilinson or
Bloch--Kato conjectures on special values of $L$-functions
necessitate
the explicit construction of elements in motivic cohomology.
Furthermore, the realizations (Betti, de~Rham,
\'etale...) of these elements can often be constructed out
from the cohomology of non-compact varieties, using the 
respective
realizations of our exact triangle
\[
\dMgm (X) \longto \Mgm (X) \longto \Mcgm (X) \longto \dMgm (X) [1] \; .
\]
This approach is clearly present e.g.\ in Harder's work on 
special values. 
Our definition may thus be seen as an attempt to give a rigorous
motivic meaning to these constructions. \\

In order to efficiently apply our new notion,
one is naturally led to look for means to identify 
the boun\-dary motive.
We develop three tools to compute $\dMgm (X)$: co-localization
(Section~\ref{3}), invariance under abstract blow-up
(Section~\ref{3a}), and analytical invariance (Section~\ref{4}).
In a sequel to this work \cite{W9}, they will be applied to
Shimura varieties, yielding in particular a motivic version of 
Pink's theorem on higher direct images of \'etale sheaves 
in the Baily--Borel compactification. \\

All our tools are based on the identification of $\dMgm (X)$ with the
motive associated to the diagram of schemes
\[
\xymatrix@R-20pt{
\emptyset \ar[r] \ar[dd] & 
X \ar[dd] \\
& \\
\doX \ar[r] & 
\oX 
\\}
\]
for a compactification $\oX$ of $X$, with $\doX := \oX - X$
(see Proposition~\ref{2C} for the precise statement).
Let us insist that the definition of $\dMgm (X)$ does not
involve such a compactification. \\

Co-localization is the motivic analogue of the dual of the
localization spectral sequence associated to a stratification
of $\doX$. In the context of Betti cohomology (say),
a good stratification $\doX = \coprod Y_m$ induces a spectral
sequence converging to the cohomology of $\doX$, and with
$E_2$-terms given by cohomology with compact support of the $Y_m$.
If one wants to express the analogue of the $E_2$-terms as
motives rather than motives with compact support, one is led to
imitate the dual spectral sequence, whose $E_2$-terms are equal
to usual cohomology of $Y_m$, but with coefficients given by the
exceptional inverse images of the coefficients on $\doX$. 
This explains our choice of notation for the motives ``with
coefficients'' occurring in Sections~\ref{3}--\ref{4}, and 
more importantly, the exact nature of their behaviour. It turns
out that in order to get a clean statement, co-localization 
(Theorem~\ref{3C}) has to be formulated more generally 
for diagrams of the type
\[
\xymatrix@R-20pt{
Y' - Y \ar[r] \ar[dd] & 
W - Y \ar[dd] \\
& \\
Y' \ar[r] & 
W 
\\}
\]
and for stratifications of $Y$,
where $Y \into Y' \into W$ are closed immersions of not necessarily
proper schemes over $k$. This diagram should be seen as modeling
cohomology of the ``motivic sheaf'' $\iYjZ$, where $i_Y: Y \into W$ and
$j: W-Y' \into W$. Due to the lack of sufficient covariance properties
of the functor $\Mcgm$, it does not seem obvious to geometrically
model cohomology with compact support of $\iastYjZ$. 
We refer to Definition~\ref{7A} for a partial re\-me\-dy to this problem. \\

Given the sheaf theoretical 
point of view, 
invariance under abstract blow-up
and analytical invariance
come as no surprise:
in the context of complex spaces (say), 
the cohomology of $\iYjZ$ can be computed after
proper base change to a second diagram as above,
provided that this base change induces an isomorphism
above the complement $W - Y'$.
Theorem~\ref{3aA} is the
motivic analogue of this invariance. 
Verdier's theory of the specialization functor shows that
the sheaf
$\iYjZ$ on $Y$ can be computed on the normal cone of $Y$ along $W$.
Since analytically isomorphic situations near $Y$ lead to the
same normal cone, this implies that $ \iYjZ$ can be computed with
respect to a second set of closed immersions $Y \into Y'_2 \into W_2$,
provided that the formal completions of $W_2$, resp. of
$Y'_2$ along $Y$ agree with those of $W$, resp. of
$Y'$ along $Y$. Theorem~\ref{4A} states that the
motivic analogue of this latter statement holds. While 
co-localization and invariance under abstract blow-up
are direct consequences of the material contained in \cite{VSF}, 
the proof of analytical 
invariance uses in addition the full force of Artin approximation. \\

In Section~\ref{6}, we generalize duality for 
bivariant cycle cohomology \cite[Thm.~IV.7.4]{VSF}, in order to establish
an important structural property of
the boundary motive $\dMgm (X)$ of a smooth scheme $X$ of pure dimension 
$n$: it is canonically isomorphic to its own dual $\dMgm (X)^*$, twisted by
$n$ and shifted by $2n-1$ (Theorem~\ref{6main}). As one 
formal consequence of this auto-duality, and of co-localization, 
we obtain a fourth tool to compute  $\dMgm (X)$, namely
localization (Section~\ref{7})
in the context of stratifications
$\doX = \coprod Y_m$, for a compactification $\oX$ of a smooth
scheme $X$ (with $\doX = \oX - X$ as above). \\

Sections~\ref{3}--\ref{6} are logically independent of each other.
Section~\ref{1} serves as basis for all that is to
follow, and Section~\ref{7} uses everything said before. 
All results from Sections~\ref{3a}, \ref{6} and \ref{7} require resolution
of singularities for the base field $k$. \\ 

This work was done while I was enjoying a \emph{d\'el\'egation
aupr\`es du CNRS}, 
and during visits to the \emph{Sonderforschungsbereich~478}
of the University of M\"unster, and to the \emph{Institut de Matem\`atica}
of the University of Barcelona. I am grateful to all three institutions.
I also wish to thank J.~Barge, C.~Denin\-ger, F.~Lemma, M.~Levine and F.~Morel for
useful discussions and comments. The remarks of the referee were appreciated;
it is thanks to her or him that Section~\ref{8} (on the case of a normal crossing
boundary $\doX$) was added.


\bigskip

%
%

\section{Notations and conventions}
\label{1}



Our main and almost only reference is the book \cite{VSF}. When citing
a result from its Chapter $n$, we shall precede the numbering
used in \loccit \ by $N$, where $N$ is the symbol representing $n$
in the Roman number
system. Example: Proposition~3.1.3 from Chapter~5 from \cite{VSF}
will be cited as \cite[Prop.~V.3.1.3]{VSF}. \\

We follow the notation of \cite{VSF}. Fix a perfect base field
$k$. Denote by $Sch/k$ the category of schemes which are
separated and of finite type
over $k$, and by $Sm/k$ the full sub-category of objects which
are smooth over $k$. Recall the definition
of the category $\SmC$ \cite[p.~190]{VSF}: 
its objects are those of $Sm/k$. Morphisms
from $Y$ to $X$ are given by the group $c(Y,X)$ of \emph{finite 
correspondences} from $Y$ to $X$, defined as the free Abelian group
on the symbols $(Z)$, where $Z$ runs through the integral closed
sub-schemes of $Y \times_k X$ which are finite over $Y$ and surjective
over a connected component of $Y$. Note for later use that the
definition of $c(Y,X)$ still makes sense when $X \in Sch/k$ is not
necessarily smooth. The category $\ShN$
of \emph{Nisnevich sheaves with transfers}
\cite[Def.~V.3.1.1]{VSF} is the category of those
contravariant additive functors from $\SmC$ to Abelian groups,
whose restriction to $Sm/k$ is a sheaf for the Nisnevich topo\-logy.
This category is Abelian \cite[Thm.~V.3.1.4]{VSF}.
Inside the derived category $D^-(\ShN)$ of complexes bounded from
above, one defines the full sub-category $\DM$ 
of \emph{effective motivic complexes} over $k$
\cite[p.~205]{VSF} as the one consisting
of objects whose cohomology sheaves are \emph{homotopy invariant}
\cite[Def.~V.3.1.10]{VSF}. Since $k$ is supposed to be perfect,
this sub-category is triangulated \cite[Prop.~V.3.1.13]{VSF}.
According to \cite[Prop.~V.3.2.3]{VSF}, the inclusion of
$\DM$ into $D^-(\ShN)$ admits a left adjoint $\RC$, which is 
induced from the functor
\[
\uC_*: \ShN \longto C^-(\ShN) 
\]
which maps $F$ to the simple complex associated to the 
\emph{singular simplicial complex} \cite[p.~207]{VSF}. 
Its $n$-th term (in homological
numbering) $\uC_n (F)$ sends
$X$ to $F(X \times_k \Delta^n)$. \\

One defines two functors $L$ and $L^c$ from $Sch/k$ to
$\ShN$ \cite[pp.~223, 224]{VSF}: the functor $L$ associates
to $X$ the Nisnevich sheaf with transfers $c(\argdot,X)$.
The functor $L^c$ maps $X$ to 
\[
Y \longmapsto z(Y,X) \; ,
\]
$z(Y,X)$ being defined as the free Abelian group
on the symbols $(Z)$, where $Z$ runs through the integral closed
sub-schemes of $Y \times_k X$ which are quasi-finite over $Y$ and 
dominant over a connected component of $Y$. One defines
the \emph{motive} $\Mgm (X)$ of $X \in Sch/k$ as $\RC (L(X))$,
and the \emph{motive with compact support} $\Mcgm (X)$ as
$\RC (L^c(X))$. \\

For certain applications, it is of interest
to enlarge the domain of the functor $L$: denote by
$\Schi$ the category of schemes which are
separated and \emph{locally} of finite type
over $k$. The functor $L$ extends, with the same definition 
of $c(\argdot,X)$ as above. This identifies
$L (X)$ with the filtered direct limit of the
$L (U)$, with $U$ running through the open sub-schemes of $X$
which are of finite type over $k$. This observation allows to use
certain results from \cite{VSF} also for the motives 
$\Mgm (X) := \RC (L(X)) \in \DM$,
with $X$ in $\Schi$. \\

We shall also use another, more geometric approach to motives, i.e.,
the one developed in \cite[V.2.1]{VSF}. There, the triangulated
category $\DeffgM$ of \emph{effective geometrical motives} over $k$
is defined. There is a canonical 
full triangulated embedding of $\DeffgM$ into $\DM$ \cite[Thm.~V.3.2.6]{VSF}, 
which maps the geometrical
motive of $X \in Sm/k$ \cite[Def.~V.2.1.1]{VSF} to $\Mgm (X)$.
Using this embedding, we consider $\Mgm (X)$ as an object of $\DeffgM$. 
Finally, the category $\DgM$ of \emph{geometrical motives} over $k$
is obtained from $\DeffgM$ by inverting the \emph{Tate motive} $\BZ(1)$
\cite[p.~192]{VSF}. 
All four categories $\DeffgM$, $\DgM$, $D^-(\ShN)$, and $\DM$
are tensor triangulated, and admit unit objects 
\cite[Prop.~V.2.1.3, Cor.~V.2.1.5, p.~206,
Thm.~V.3.2.6]{VSF}. These tensor structures are such that
for all $X, Y \in Sm/k$, one has
\[
\Mgm (X) \otimes \Mgm (Y) = \Mgm (X \times_k Y)
\]
in $\DeffgM$, and
\[
L(X) \otimes L(Y) = L(X \times_k Y)
\]
in $D^-(\ShN)$. The unit object of $\DeffgM$ is $\Mgm (\Spec k)$, and
that of $D^-(\ShN)$ is $L (\Spec k)$.
Both of them are denoted by $\BZ (0)$.
For $M \in \DgM$ and $n \in \BZ$, write $M(n)$ for the tensor product
$M \otimes \BZ (n)$.
The three
functors $\DeffgM \into \DM$, $\DeffgM \to \DgM$, and
$\RC: D^-(\ShN) \to \DM$ are compatible with the tensor structure.
(By contrast, the embedding of $\DM$ into $D^-(\ShN)$ is not;
see \cite[Remark on p.~206]{VSF}.) 
According to \cite[Thm.~V.4.3.1]{VSF},
the functor $\DeffgM \to \DgM$ is a full triangulated
embedding if $k$ admits resolution of singularities.

\begin{Con}
Whenever we speak about resolution of singularities,
it will be taken in the sense of 
\cite[Def.~IV.3.4]{VSF}. 
\end{Con}

\begin{Con} \label{1A}
We shall use the same symbol for $\Mgm (X) \in \DM$ and for
its canonical representative $\uC_* (L(X))$ in $C^- (\ShN)$;
similarly  for $\Mcgm (X)$. Whenever we speak about cones
of morphisms between motives, we mean the class of the cone
of the morphism between the canonical representatives.
For a commutative diagram 
\[
\xymatrix@R-20pt{
& X' \ar[r] \ar[dd] & 
X_1 \ar[dd] \\
\FX \quad = & & \\
& X_2 \ar[r] & 
X'' 
\\}
\]
in $\Schi$, define its motive $\Mgm (\FX) \in \DM$ as $\RC$
applied to the simple complex $s L(\FX)$ associated to $L(\FX)$,
i.e., to the complex 
\[
L(X') \longto L(X_1) \oplus L(X_2) \longto L(X'') \; ,
\]
which we normalize by assigning degree zero to the component 
$L(X'')$. A si\-mi\-lar
construction is possible for commutative
diagrams in $\Schi$ of ``dimension'' greater than two, provided
that there are not more than two schemes on any of the lines
in the diagram.
\end{Con}


\bigskip

%
%

\section{Definition of the boundary motive}
\label{2}



Let $X \in Sch/k$. Note that the inclusion
$c(\argdot,X) \into z(\argdot,X)$ induces a monomorphism
\[
\iota_X: L(X) \longinto L^c(X) \; .
\]

\begin{Def} \label{2A}
The \emph{boundary motive} of $X$ is defined as
\[
\dMgm (X) := \RC (\coker \iota_X) [-1] \; .
\]
\end{Def}

Note that there is a canonical quasi-isomorphism
\[
\cone \left( \Mgm (X) \longto \Mcgm (X) \right) \longto \dMgm (X) \; ,
\]
where the cone is to be understood as in Convention~\ref{1A}.
We have:

\begin{Prop} \label{2B}
There is an exact triangle
\[
\dMgm (X) \longto \Mgm (X) \longto \Mcgm (X) \longto \dMgm (X) [1] 
\]
in $\DM$. 
\end{Prop}

\begin{Cor} \label{2Ba}
Assume that $k$ admits resolution of
singularities.
Then the boundary motive $\dMgm (X)$ belongs to $\DeffgM$.
\end{Cor}

\begin{Proof}
This follows from Proposition~\ref{2B}, the fact that the embedding
of $\DeffgM$ into $\DM$ is triangulated, and 
\cite[Cor.~V.4.1.4 and Cor.~V.4.1.6]{VSF}.
\end{Proof}

The definition of the boundary motive does not involve a
compactification of $X$. However, the tools to
compute $\dMgm (X)$ which we shall develop 
in the sequel are based on the following:

\begin{Prop} \label{2C}
Let $\oX$ be a compactification of $X \in Sch/k$, and define $\doX$ as
the complement $\oX - X$, equipped 
with the reduced scheme structure. There is a canonical morphism
\[
\cone \left( \Mgm \left( X \coprod \doX \right) 
\to \Mgm (\oX) \right) \longto \dMgm (X) 
\]
in $C^-(\ShN)$, 
which becomes an isomorphism in $\DM$ if $k$ admits resolution of
singularities.
\end{Prop}   

\begin{Proof}
Consider the exact sequences 
\[
0 \longto L(X) \longto L^c(X) \longto \coker \iota_X \longto 0 \; ,
\]
and 
\[ 
0 \longto L^c(\doX) \longto L^c(\oX) \longto L^c(X) \; .
\] 
Observe that since $\doX$ and $\oX$ are proper, we have
$L(\doX) = L^c(\doX)$ and $L(\oX) = L^c(\oX)$. The monomorphism
$L(X) \into L^c(X)$ factors through $L(\oX)$. Hence the exact
sequences induce a monomorphism between the quotient
$L(\oX) / (L(X \coprod \doX)$ and $\coker \iota_X$, whose
cokernel is identical to that of the restriction
$L(\oX) = L^c(\oX) \to L^c(X)$.
According to \cite[Prop.~V.4.1.5]{VSF}, this latter cokernel has trivial image
under $\RC$ if $k$ admits resolution of singularities.
\end{Proof}


\bigskip

%
%

\section{Co-localization}
\label{3}



Consider the geometric situation of Proposition~\ref{2C}: let $X \in Sch/k$,
choose a compactification $\oX$ of $X$, and define $\doX := \oX - X$.
In this section, we develop the motivic analogue of the dual of the
localization spectral sequence associated to a stratification of $\doX$. 
It turns out to be useful to consider a more general geometric situation:
fix closed immersions $Y \into Y' \into W$ in $\Schi$. Write
$j$ for the open immersion of $W-Y'$, and $i_Y$ for the closed
immersion of $Y$ into $W$. Denote by $\FY$ the commutative diagram
\[
\xymatrix@R-20pt{
Y' - Y \ar[r] \ar[dd] & 
W - Y \ar[dd] \\
& \\
Y' \ar[r] & 
W 
\\}
\]

\begin{Def} \label{3A}
The \emph{motive of $Y$ with
coefficients in $\iYjZ$} is defined as
\[
\Mgm (Y, \iYjZ) := \Mgm(\FY) 
\]
(see Convention~\ref{1A}).
\end{Def}
  
\begin{Rem} \label{3Aa}
(a)~Note that $\Mgm (Y, \iYjZ)$ lies in $\DeffgM$ if $W$ and $Y'$
are in $Sm/k$. \\[0.1cm]
(b)~Assume that $k$ admits resolution of singularities.
Then $\Mgm (Y, \iYjZ)$ lies in $\DeffgM$ if $W$
is in $Sch/k$
\cite[Cor.~V.4.1.4]{VSF}.
\end{Rem}

\begin{Rem} \label{3B}
If $Y=Y'$, i.e., if $i_Y$ is complementary to $j$, 
then $\FY$ acquires the form 
\[
\xymatrix@R-20pt{
\emptyset \ar[r] \ar[dd] & 
W - Y \ar[dd] \\
& \\
Y \ar[r] & 
W 
\\}
\]
and we have 
\[
\Mgm (Y, \iYjZ) = 
\cone \left( \Mgm \left( (W - Y) \coprod Y \right) 
\to \Mgm (W) \right) [1] \; .   
\]
If in addition $W$ is proper, then Proposition~\ref{2C} relates
$\Mgm (Y, \iYjZ)$ to $\dMgm (W-Y)[1]$. 
\end{Rem}

Now assume given a filtration
\[
\emptyset = \FF_{-1} Y \subset \FF_0 Y \subset \ldots \subset \FF_d Y = Y
\]
of $Y$ by closed sub-schemes. It induces a stratification of $Y$ by
locally closed sub-schemes $Y_m := \FF_m Y - \FF_{m-1} Y$, for
$m = 0, \ldots, d$. Define $W^m$ as the complement of $\FF_{m-1} Y$ in
$W$. This gives a descending partial filtration of $W$ by open sub-schemes.
Note in particular that we have $W^0 = W$ and $W^{d+1} = W - Y$.
Write $i_{Y_m}$ for the closed immersion of $Y_m$ into $W^m$.
By abuse of notation, we use the letter $j$ to denote also the open
immersions of $W-Y'$ into $W^m$.

\begin{Thm}[Co-localization] \label{3C}
There is a cano\-ni\-cal chain of morphisms
\[
M^{d+1} = 0 \stackrel{\gamma^d}{\longto} M^d  
\stackrel{\gamma^{d-1}}{\longto} M^{d-1}
\stackrel{\gamma^{d-2}}{\longto} \ldots
\stackrel{\gamma^0}{\longto} M^0 = \Mgm (Y, \iYjZ)
\]
in $C^- (\ShN)$. For each 
$m \in \{ 0, \ldots, d\}$, there is a canonical isomorphism
\[
\cone \gamma^m \cong \Mgm (Y_m, \iYmjZ) [-1]
\]
and hence, a canonical exact triangle
\[
\Mgm (Y_m, \iYmjZ) [-1] \longto M^{m+1} 
\stackrel{\gamma^m}{\longto} M^m \longto \Mgm (Y_m, \iYmjZ)  
\]
in $\DM$.
In particular, all the $M^m$ represent objects in $\DM$.
If $k$ admits resolution of singularities, then all 
the $M^m$ represent objects in $\DeffgM$.
\end{Thm}
       
\begin{Proof}
Consider the induced filtration 
$Y'^m := Y' - \FF_{m-1} Y = W^m \cap Y'$, and define
$M^m := \Mgm (\FY^m)$, where $\FY^m$ is the commutative diagram
\[
\xymatrix@R-20pt{
Y' - Y \ar[r] \ar[dd] & 
W - Y \ar[dd] \\
& \\
Y'^m \ar[r] & 
W^m 
\\}
\]
Now note that $\Mgm (Y_m, \iYmjZ)$ is the motive associated to the
diagram
\[
\xymatrix@R-20pt{
Y'^{m+1} \ar[r] \ar[dd] & 
W^{m+1} \ar[dd] \\
& \\
Y'^m \ar[r] & 
W^m 
\\}
\]
\end{Proof}

\begin{Cor} \label{3D}
Assume that $k$ admits resolution of singularities.
In the above situation, assume that $Y = \doX := \oX - X$,
with $X \in Sch/k$, and $\oX$ 
a compactification of $X$.
Write $(\doX)_m:= Y_m$.  
Then there is a cano\-ni\-cal chain of morphisms
\[
M^{d+1} = 0 \stackrel{\gamma^d}{\longto} M^d  
\stackrel{\gamma^{d-1}}{\longto} M^{d-1}
\stackrel{\gamma^{d-2}}{\longto} \ldots
\stackrel{\gamma^0}{\longto} M^0 = \dMgm (X)
\]
in $C^- (\ShN)$. 
For each 
$m \in \{ 0, \ldots, d\}$, there is a canonical exact triangle
\[
\Mgm ((\doX)_m, \iXmjZ) [-1] \to M^{m+1} 
\stackrel{\gamma^m}{\longto} M^m \to \Mgm ((\doX)_m, \iXmjZ)
\]
in $\DgM$. In particular, all the $M^m$ are in $\DeffgM$.
\end{Cor}

\begin{Proof}
This follows from Proposition~\ref{2C} and
Theorem~\ref{3C}.
\end{Proof}

Of particular interest is the case when $X \in Sm/k$,
and $\oX$ is a smooth compactification of $X$, such that
$\doX$ is a normal crossing divisor. Then $d < \dim X$, and
$(\doX)_m$ equals the geometric locus of points lying on
exactly $d+1-m$ local irreducible components of
$\doX$. Since analytical invariance will permit us to
give a good description of the $\Mgm ((\doX)_m, \iXmjZ)$, 
we postpone the discussion of this case until Section~\ref{8}.


\bigskip

%
%

\section{Invariance under abstract blow-up}

\label{3a}



Fix a proper morphism $\pi: W_1 \to W_2$, and closed immersions 
$Y_2 \into Y'_2 \into W_2$ in $\Schi$.
Write $Y_1 \into Y'_1 \into W_1$ for the base change of these
immersions, $j_m$ for the open immersion of $W_m-Y'_m$, 
and $i_{Y_m}$ for the closed
immersion of $Y_m$ into $W_m$~:
\[
\xymatrix@R-20pt{
Y_1 \ar@{^{ (}->}[r]^-{i_{Y_1}} \ar[dd] &
W_1 \ar[dd]_{\pi} &
W_1 - Y'_1 \ar@{_{ (}->}[l]_-{j_1} \ar[dd]^{\pi} \\
&& \\
Y_2 \ar@{^{ (}->}[r]^-{i_{Y_2}} &
W_2 &
W_2 - Y'_2 \ar@{_{ (}->}[l]_-{j_2}
\\}
\]

\begin{Thm}[Invariance 
under abstract blow-up]  \label{3aA} 
Assume that
\[
\pi: W_1-Y'_1 \longto W_2-Y'_2
\]
is an isomorphism.
If $k$ admits resolution of singularities,
then the map
\[
\pi: \Mgm (Y_1, i_{Y_1}^! j_{1 !} \, \BZ) \longto 
\Mgm (Y_2, i_{Y_2}^! j_{2 !} \, \BZ)
\]
is an isomorphism.
\end{Thm}

Recall that $\Mgm (Y_m, i_{Y_m}^! j_{m !} \, \BZ)$
is associated to the diagram 
\[
\xymatrix@R-20pt{
Y'_m - Y_m \ar[r] \ar[dd] & 
W_m - Y_m \ar[dd] \\
& \\
Y'_m \ar[r] & 
W_m 
\\}
\]
Since the two rows are of the same nature, Theorem~\ref{3aA} is
a formal consequence of the following:

\begin{Thm} \label{3aB} 
Consider closed immersions $Y'_1 \into W_1$ and
$Y'_2 \into W_2$ in $\Schi$, and a proper morphism
$\pi: W_1 \to W_2$ identifying $Y'_1$ with the fibre
product $W_1 \times_{W_2} Y'_2$, and inducing an isomorphism
from $W_1-Y'_1$ to $W_2-Y'_2$. 
If $k$ admits resolution of singularities,
then the monomorphism
\[
\pi: L(W_1) / L(Y'_1) \longto 
L(W_2) / L(Y'_2) 
\]
in $\ShN$ induces an isomorphism $\RC (\pi)$ in $\DM$.
\end{Thm}

\begin{Proof}
(See \cite[Prop.~V.4.1.3]{VSF} and its proof.)
The sequence
\[
0 \longto L(Y'_1) \longto L(W_1) \oplus L(Y'_2) \longto L(W_2)
\]
is exact, and the quotient $L(W_2) / (L(W_1) + L(Y'_2))$
satisfies the condition of \cite[Thm.~V.4.1.2]{VSF}.
\end{Proof}


\bigskip

%
%

\section{Analytical invariance}
\label{4}



Fix closed immersions $Y \into Y'_1 \into W_1$ and
$Y \into Y'_2 \into W_2$ in $\Schi$.
Write
$j_m$ for the open immersion of $W_m-Y'_m$, and $i_{Y,m}$ for the closed
immersion of $Y$ into $W_m$. 
The aim of this section is to prove the following:

\begin{Thm}[Analytical invariance] \label{4A}
Assume gi\-ven an isomorphism
\[
f: (W_1)_Y \isoto (W_2)_Y
\]
of formal completions along $Y$
inducing an isomorphism $(Y'_1)_Y \cong (Y'_2)_Y$, 
and compatible with the immersions $i_{Y,m}$ of $Y$.
Then $f$ induces an isomorphism 
\[
\Mgm (Y, i_{Y,1}^! j_{1 !} \, \BZ) \isoto 
\Mgm (Y, i_{Y,2}^! j_{2 !} \, \BZ) 
\]
in $\DM$.
\end{Thm}

\begin{Rem} \label{4Ba}
Using Proposition~\ref{2C} together with Remark~\ref{3B}, we deduce
the following statement
from Theorem~\ref{4A}, assuming that $k$ admits resolution of 
singula\-ri\-ties: let $\oX_m$ be
a compactification of $X_m \in Sch/k$, 
$m = 1,2$, and set $\doX_m := \oX_m - X_m$
(with the reduced scheme structure). Assume that there is an
isomorphism $\doX_1 \cong \doX_2$, which can be extended to an 
isomorphism between the formal completions of $\oX_m$ along $\doX_m$.
Then $\dMgm(X_1)$ and $\dMgm(X_2)$ are isomorphic. \\

Note however that
in practice, it may not always 
be possible to identify the formal completion
of $\doX$ along a given compactification $\oX$ of a scheme $X$.
Actually, one might control the formal
completion of an abstract blow-up of each stratum belonging to a 
stratification of 
$\doX$.
In order to compute $\dMgm (X)$ in such a situation, one first
applies co-localization with respect to the stratification,
then uses invariance under abstract blow-up for each stratum,
and finally  
analytical invariance as stated in the above generality.  
\end{Rem}

The main technical ingredient of the proof of Theorem~\ref{4A} is
the following consequence of Artin approximation:

\begin{Thm}[Artin] \label{4C}
Let $S$ be the spectrum of a field or of an excellent Dedekind
domain, $W_1$ and $W_2$ two $S$-schemes which are locally of finite type, 
$Y_m$ closed sub-schemes of $W_m$,
$y_m \in Y_m$, and $\Fa_m$ the ideal of $Y_m$ in the Henselization
$\CO_{W_m,y_m}^h$ of the local ring $\CO_{W_m,y_m}$. Denote by
$\widehat{\CO}_{W_m,y_m}$ the completion of $\CO_{W_m,y_m}^h$ with respect
to $\Fa_m$, for $m = 1,2$. \\[0.1cm]
(a) If $\widehat{\CO}_{W_1,y_1} \cong \widehat{\CO}_{W_2,y_2}$ over $\CO_S$,
then $y_1$ and $y_2$ have a common Nisnevich neighbourhood: there exists
an $S$-scheme $W'$, a point $y' \in W'$, and \'etale morphisms
$W' \to W_m$ mapping $y'$ to $y_m$, for $m= 1,2$, which identify the
residue fields $\kappa(y_1) \cong \kappa(y') \cong \kappa(y_2)$. \\[0.1cm] 
(b) Assume that in addition we are given an isomorphism
$Y_1 \cong Y_2$. Assume that
the isomorphism in (a) maps the completed ideal
$\hat{\Fa}_1$ isomorphically to $\hat{\Fa}_2$, and that the induced
isomorphism $\CO_{Y_1,y_1}^h \cong \CO_{Y_2,y_2}^h$ is compatible with
the given isomorphism
$Y_1 \cong Y_2$. Then the Nisnevich neighbourhood $W'$
in (a) can be chosen such that in addition 
\[
Y' := W' \times_{W_1} Y_1 = W' \times_{W_2} Y_2
\]
as sub-schemes of $W'$, and the induced \'etale morphisms $Y' \to Y_m$
are compa\-tible with the isomorphism $Y_1 \cong Y_2$.  
\end{Thm}

\begin{Proof}
This is a variant of \cite[Cor.~(2.6)]{A}. In fact, the results stated
in Section~2 of \loccit \ are the translations of the main results
of Section~1 only in the case when the ideal of definition
is the maximal ideal of the point in question. In order to deduce
the variant from \cite[Thm.~(1.12)]{A}, one faithfully imitates
the proof of \cite[Cor.~(2.6)]{A}.  
\end{Proof}

\begin{Cor} \label{4D}
With $S$, $W_m$ and $Y_m$ as in Theorem~\ref{4C}, assume given 
an isomorphism $Y_1 \cong Y_2$, which extends to an isomorphism
\[
f: (W_1)_{Y_1} \isoto (W_2)_{Y_2}
\]
of formal completions. Then there are Nisnevich coverings $\FW_m$
of $W_m$ of the form
\[
\FW_1 = \{ W'_i \tei i \in I \} \coprod \{ W_1 - Y_1 \} \quad , \quad
\FW_2 = \{ W'_i \tei i \in I \} \coprod \{ W_2 - Y_2 \}
\]
(with the same $W'_i$~!) such that for any $i \in I$, one has
\[
Y'_i := W'_i \times_{W_1} Y_1 = W'_i \times_{W_2} Y_2
\]
as sub-schemes of $W'_i$, and the induced \'etale morphisms $Y'_i \to Y_m$
are compa\-tible with the isomorphism $Y_1 \cong Y_2$. 
\end{Cor}

We turn to the proof of Theorem~\ref{4A}.
Recall that $\Mgm (Y, i_{Y,m}^! j_{m !} \, \BZ)$
is associated to the diagram 
\[
\xymatrix@R-20pt{
Y'_m - Y \ar[r] \ar[dd] & 
W_m - Y \ar[dd] \\
& \\
Y'_m \ar[r] & 
W_m 
\\}
\]
Since the two columns are of the same nature, Theorem~\ref{4A} is
a formal consequence of parts~(a) and (b) of the following:

\begin{Thm} \label{4E} 
Consider closed immersions $Y \into W_1$ and
$Y \into W_2$ in $\Schi$, which extend to
an isomorphism
\[
f: (W_1)_Y \isoto (W_2)_Y
\]
of formal completions along $Y$. \\[0.1cm]
(a) There is an isomorphism
\[
L(W_1) / L(W_1 - Y) \isoto 
L(W_2) / L(W_2 - Y) 
\]
of Nisnevich sheaves with transfers, and depending only on $f$. \\[0.1cm]
(b)~The isomorphism $L(W_1) / L(W_1 - Y) \isoto L(W_2) / L(W_2 - Y)$
is compatible with restriction of the $W_m$ to sub-schemes
containing $Y$. \\[0.1cm]
(c)~The isomorphism $L(W_1) / L(W_1 - Y) \isoto L(W_2) / L(W_2 - Y)$
is  compatible with restriction of $Y$.
\end{Thm}
\forget{
\begin{Rem}
The functoriality statements
(b) and (c) of Theorem~\ref{4E} imply that analytical invariance
and co-localization are well compatible. More precisely, 
in the situation of Theorem~\ref{4A}, assume in addition that 
a filtration
\[
\emptyset = \FF_{-1} Y \subset \FF_0 Y \subset \ldots \subset \FF_d Y = Y
\]
of $Y$ by closed sub-schemes is given. According to Variant~\ref{3D},
the motives $\Mgm(Y_m, \iYmjZ)$ can be computed from the diagrams
$\FY^{(3)}_m$ associated to this filtration, and formed with
respect to the ambient scheme $W_m$, $m = 1,2$.
Then according to Theorem~\ref{4E}, the isomorphism from Theorem~\ref{4A}
can be realized by an isomorphism 
\[
sL(\FY^{(3)}_1) \isoto sL(\FY^{(3)}_2)
\]
depending only on $f$.
\end{Rem}
}

\begin{Proof}
Choose Nisnevich coverings 
\[
\FW_m = \{ W'_i \tei i \in I \} \coprod \{ W_m - Y \}
\] 
as in 
Corollary~\ref{4D}.
Set
\[
W' := \coprod_{i \in I} W'_i \; ,
\] 
and write $\alpha_m$ for the coproduct of the \'etale morphisms from
the $W'_i$ to $W_m$, for $m = 1,2$.
By Corollary~\ref{4D}, we have
\[
Y' := W' \times_{W_1} Y = W' \times_{W_2} Y \; ,
\]
and $\alpha_1$ and $\alpha_2$ coincide on $Y'$.
Using \cite[Prop.~V.3.1.3]{VSF}, we see that we have exact sequences
\[
\frac{L(W' \times_{W_1} W')}{L((W'-Y') \times_{W_1} (W'-Y'))}
\stackrel{pr^1_* - pr^2_*}{\longto}
\frac{L(W')}{L(W'-Y')} 
\stackrel{\alpha_{1 *}}{\longto} 
\frac{L(W_1)}{L(W_1-Y)} \longto 0 \; ,
\]
\[
\frac{L(W' \times_{W_2} W')}{L((W'-Y') \times_{W_2} (W'-Y'))}
\longto
\frac{L(W')}{L(W'-Y')} 
\stackrel{\alpha_{2 *}}{\longto} 
\frac{L(W_2)}{L(W_2-Y)} \longto 0 \; ,
\]
of Nisnevich sheaves with transfers. Let us show that the map
$\alpha_{2 *}$ is zero on the
image of $L(W' \times_{W_1} W') / L((W'-Y') \times_{W_1} (W'-Y'))$.
We imitate the proof of \cite[Prop.~II.4.3.9]{VSF}. It is
sufficient to show the following claim:
\begin{itemize}
\item[$(*)$] For any local Henselian scheme $S$ which is smooth
over $k$, the composition
\[
c(S,W' \times_{W_1} W') 
\stackrel{pr^1_* - pr^2_*}{\longto}
c(S,W') \stackrel{\alpha_{2 *}}{\longto} 
\frac{c(S,W_2)}{c(S,W_2-Y)}
\]
is trivial. 
\end{itemize}
Note that the presheaves $U \mapsto c(U,T)$, for
$T \in \Schi$, can be extended in an obvious way to the category
of smooth $k$-schemes which are not necessarily of finite type.
For the proof of $(*)$, we shall repeatedly apply the following
principle, valid since $S$ is Henselian: for any $T \in \Schi$,
the support of any element of $c(S,T)$ is a disjoint union of
local Henselian schemes. This principles reduces us to consider
only cycles in $c(S,W' \times_{W_1} W')$ of the form $(Z)$,
where $Z$ is a local Henselian sub-scheme of $S \times_k (W' \times_{W_1} W')$.
Without loss of generality, we may assume that the closed point
of $Z$ lies over $Y'$, hence over $Y$.
Write $pr^l_* (Z) = n^l \cdot (Z^l)$, with local Henselian sub-schemes
$Z^l$ of $S \times_k W'$, for $l = 1,2$. We have 
$\alpha_{1 *} pr^1_* = \alpha_{1 *} pr^2_*$, hence the 
$\alpha_{1 *} (Z^l)$ are multiples of $(Z_1)$, for \emph{one} local
Henselian sub-scheme $(Z_1)$ of $S \times_k W_1$. In order to show the analogous
statement for the $\alpha_{2 *} (Z^l)$, note first that the closed point
$y$ of $Z_1$ belongs to $S \times_k Y$. The support of
\[
\alpha_{2 *} (pr^1_* - pr^2_*) (Z) = 
\alpha_{2 *} (n^1 \cdot (Z^1) - n^2 \cdot (Z^2))
\]
is a disjoint union of local Henselian schemes, parametrized by their closed
points. But since $\alpha_1$ and $\alpha_2$ coincide on $Y'$,
this support must be local, and we have indeed
\[
\alpha_{2 *} (pr^1_* - pr^2_*) (Z) = 
r \cdot (Z_2) \; ,
\]
for a local Henselian sub-scheme $(Z_2)$ of $S \times_k W_2$,
whose closed point is $y$. In order to show that $r = 0$,
consider the commutative diagram
\[
\xymatrix@R-20pt{
& c(S,W') \ar[dl]_-{\alpha_{1 *}} \ar[dr]^-{\alpha_{2 *}} & \\
c(S,W_1) \ar[dr]_-{w_{1 *}} & & c(S,W_2) \ar[dl]^-{w_{2 *}} \\
& c(S,\Spec k) & 
\\}
\]
$w_m$ denoting the structure morphism of $W_m$. 
On the one hand,
\[
w_{2 *} \alpha_{2 *} (pr^1_* - pr^2_*) (Z) = 
w_{1 *} \alpha_{1 *} (pr^1_* - pr^2_*) (Z) = w_{1 *}(0) = 0 \; . 
\]
On the other hand, $w_{2 *} (Z_2)$ is non-zero
since $Z_2$ is finite over $S$. Hence
$r$ must indeed be zero, and thus
\[
\alpha_{2 *} (pr^1_* - pr^2_*) (Z) = 0 \; .
\]
This shows that $\alpha_{2 *}$ is zero on the image
of $pr^1_* - pr^2_*$.
By symmetry, we see that the identity on $L(W') / L(W'-Y')$
factors to give an isomorphism
\[
L(W_1) / L(W_1 - Y) \isoto 
L(W_2) / L(W_2 - Y) \; . 
\]
In order to prove that it does not depend on the choice of the
Nisnevich coverings $\FW_m$ as in Corollary~\ref{4D}, 
use the fact that the system of such coverings is filtering. 
\end{Proof}


\bigskip

%
%

\section{Auto-duality}
\label{6}



Throughout this section, we assume that $k$ admits resolution of singularities.
Under this assumption, $\DgM$ is a rigid tensor triangulated category
\cite[Thm.~V.4.3.7~1.\ and 2.]{VSF}. 
In particular, there exists an internal $Hom$
functor 
\[
\uHom: \DgM \times \DgM \longto \DgM \; .
\]
Writing $M^* := \uHom (M, \BZ (0))$,
we thus have $M = (M^*)^*$ for all $M \in \DgM$. \\

Now fix $X \in Sm/k$, and assume that $X$ is of pure dimension $n$.
According to \cite[Thm.~V.4.3.7~3.]{VSF}, there is a canonical 
isomorphism
\[
\mu_X: \Mcgm (X) \isoto \Mgm (X)^* (n)[2n] \; ,
\]
hence by duality, a canonical isomorphism
\[
\nu_X := \mu_X^* (n)[2n]: \Mgm (X) \isoto \Mcgm (X)^* (n)[2n] \; .
\]
The aim of this section is to prove the following:

\begin{Thm}[Auto-duality] \label{6main}
There exists a canonical isomorphism
\[
\eta_X: \dMgm (X) \isoto \dMgm (X)^* (n)[2n-1] \; .
\]
It fits into a morphism of exact triangles
\[
\xymatrix@R-20pt{
\dMgm (X) \ar[r]^-{\eta_X} \ar[dd]_{\alpha} &
\dMgm (X)^*(n)[2n-1] \ar[dd]^{\gamma^*(n)[2n]} \\
& \\
\Mgm (X) \ar[r]^-{\nu_X} \ar[dd]_{\beta} &
\Mcgm (X)^*(n)[2n] \ar[dd]^{\beta^*(n)[2n]} \\
& \\
\Mcgm (X) \ar[r]^-{\mu_X} \ar[dd]_{\gamma} &
\Mgm (X)^*(n)[2n] \ar[dd]^{\alpha^*(n)[2n]} \\
& \\
\dMgm (X)[1] \ar[r]^-{\eta_X[1]} &
\dMgm (X)^*(n)[2n]
\\}
\]
Furthermore, it is itself auto-dual in the sense that the equality
\[
\eta_X = \eta_X^* (n)[2n-1]
\]
holds.
\end{Thm}

First observe that the existence of \emph{some} isomorphism
$\eta_X$ fitting into the above diagram of exact triangles is a consequence
of the axioms of triangulated categories, and the 
commutativity of 
\[
\xymatrix@R-20pt{
\Mgm (X) \ar[r]^-{\nu_X} \ar[dd]_{\beta} &
\Mcgm (X)^*(n)[2n] \ar[dd]^{\beta^*(n)[2n]} \\
& \\
\Mcgm (X) \ar[r]^-{\mu_X} &
\Mgm (X)^*(n)[2n] 
\\}
\]
Thus, the point of auto-duality is that $\eta_X$ can be
defined canonically, and that this definition is itself
auto-dual. \\

For the proof of Theorem~\ref{6main}, 
observe that by adjunction, the construction of $\eta_X$
is equi\-va\-lent to the construction of a pairing
\[
(\argdot,\argdot): \dMgm (X) \otimes \dMgm (X) [1] \longto \BZ (n)[2n] \; .
\]
We are thus led to investigate morphisms in $\DeffgM$ whose target
is $\BZ (n)[2n]$. The statement we are aiming at is Theorem~\ref{6L}. 
It is a ge\-ne\-ralization of \cite[Cor.~V.4.2.5]{VSF}. In the above
geometric context, it implies that certain 
codimension $n$-cycles on the self product $\oX \times_k \oX$, where
$\oX$ is a smooth compactification of $X$, 
yield morphisms of the type of $(\argdot,\argdot)$. 
In order to prepare Theorem~\ref{6L},
we need to prove a variant of duality for bivariant cycle cohomology
(Theorem~\ref{6dual}). Let us start by recalling the definition
of certain variants
of the Nisnevich sheaves with transfers $L^c (V) = z(\argdot,V)$
(see \cite[p.~228]{VSF}):

\begin{Def} \label{6A}
Let $V \in Sch/k$, $W \in Sm/k$, and $r \ge 0$. \\[0.1cm]
(a)~The Nisnevich sheaf with
transfers $\zeq (V,r)$ associates to
$U \in Sm/k$ the free Abelian group
on the symbols $(Z)$, where $Z$ runs through the integral closed
sub-schemes of $U \times_k V$ which are 
equidimensional of relative dimension $r$ over $U$ and 
dominant over a connected component of $U$. \\[0.1cm]
(b)~The Nisnevich sheaf with
transfers $\zeq (W,V,r)$ 
maps $U \in Sm/k$ to $\zeq(V,r)(U \times_k W)$.
\end{Def}

Note that the sheaves $\zeq (W,V,r)$ are contravariant in the first
variable.
Recall \cite[Cor.~V.4.1.8]{VSF}:

\begin{Prop} \label{6B}
The object $\BZ(n) [2n]$ of $\DeffgM \subset \DM$ is re\-presented
by the complex
\[
\Mcgm (\BA^n) = \uC_* (z (\argdot, \BA^n)) = \uC_* (\zeq (\BA^n, 0)) \; .
\]
\end{Prop}

More generally, if $W \in Sm/k$, then
the complex
\[
\uC_* (\zeq (W, \BA^n, 0)) 
\] 
represents the functor
\[
\uHom_{\DM} (\Mgm (W), \BZ(n)[2n]) 
\]
on $\DM$ \cite[Cor.~V.4.2.7]{VSF}.
We need a variant of this statement: 

\begin{Def} \label{6CA}
Let $V \in Sch/k$, $W \in Sm/k$, and $r \ge 0$. 
Assume further that $Y \subset W$ is a closed sub-scheme. \\[0.1cm] 
(a)~Define
\[
\Mgm (W / Y):= \uC_* (L(Y \to W)) \; ,
\]
where $L (Y \to W)$ the complex given by
$L (Y)$ in degree $-1$ and $L(W)$ in degree
zero, the differential being induced by the immersion of $Y$. \\[0.1cm]
(b)~Assume in addition that arbitrary intersections of
the components $Y_j$ of $Y$ are smooth. Define 
\[
\zeq (Y^\bullet \to W, V, r)
\]
as the complex of Nisnevich sheaves with
transfers whose zeroth component is
$\zeq (W,V,r)$, and whose $m$-th component, for $m \ge 1$,
is the direct sum of the 
$\zeq (Y_J,\BA^n,0)$, for all $m$-fold intersections $Y_J$ of the
$Y_j$. 
The differentials are induced by contravariance of the 
sheaves $\zeq (\argdot,V,r)$.
\end{Def}

Of course, $\Mgm (W / Y)$ is the same thing as $\Mgm$ evaluated at the
diagram
\[
Y \longto W
\]
(see Convention~\ref{1A}).

\begin{Prop} \label{6CB}
Let $W \in Sm/k$.
Assume further that $Y \subset W$ is a closed sub-scheme
such that arbitrary intersections of
the components of $Y$ are smooth. Then the complex
\[
\uC_* (\zeq (Y^\bullet \to W, \BA^n, 0)) 
\]
represents the functor
\[
\uHom_{\DM} (\Mgm (W / Y), \BZ(n)[2n]) 
\]
on $\DM$. This identification is compatible with passage from
the pair $Y \subset W$ to $Y' \subset U$, for open
sub-schemes $U$ of $W$, and closed sub-schemes $Y'$ of $Y \cap U$
such that arbitrary intersections of
the components $Y_j'$ of $Y'$ are smooth.
\end{Prop}  

\begin{Con} \label{6conv}
By over-simplification of language,
we shall refer to the last compatibility statement in Proposition~\ref{6CB}
as ``compatibility with restriction of $W$ and $Y$''.
\end{Con} 

\begin{Proofof}{Proposition~\ref{6CB}}
First, observe that  
the canonical morphism $\Mgm (Y^\bullet) \to \Mgm (Y)$
is an isomorphism. This follows from
induction on the number of components $Y_j$. The induction
step is provided by \cite[Prop.~V.4.1.3]{VSF}.
Similarly, $\Mgm (Y^\bullet \to W) \to \Mgm (W / Y)$ is
an isomorphism.

As in
the proof of \cite[Prop.~V.4.2.8]{VSF}, one has a canonical morphism of
complexes $can$ from 
$\uC_* (\zeq (Y^\bullet \to W, \BA^n, 0)) =
\uC_* (p_*p^*(L^c(\BA^n)))$ to 
$Rp_* (p^*\uC_* (L^c(\BA^n)))$,
where $p$ denotes the structure morphism of the diagram $Y^\bullet \to W$.
By an obvious generalization of the last part of \cite[Prop.~V.3.2.8]{VSF},
one has
\[
Rp_* (p^*\uC_* (L^c(\BA^n))) = \uHom (\Mgm (Y^\bullet \to W), \BZ (n) [2n])
\]
in $\DM$. 
To check that $can$ is an isomorphism, one uses the spectral sequences
on both its source and target,
associated to the stupid filtration of $Y^\bullet \to W$.
\cite[Cor.~V.4.2.7]{VSF} shows that $can$ is an isomorphism on the
$E_1$-terms of this spectral sequence.
\end{Proofof}

\begin{Rem}
In the situation of Proposition~\ref{6CB}, assume in addition that
$W$ is of dimension at most $n$.
By \cite[Cor.~V.4.3.6]{VSF},   
the object $\Mgm (W / Y)^*(n)[2n]$ of $\DgM$ 
belongs to $\DeffgM$, and
its image under the embedding into $\DM$ equals
$\uHom_{\DM} (\Mgm (W / Y), \BZ(n)[2n])$. 
It is thus represented by the complex
$\uC_* (\zeq (Y^\bullet \to W, \BA^n, 0))$.
\end{Rem}

Given an object $F$ of $\ShN$, denote 
by $\uh^l (F)$ the $l$-th cohomology object of the complex
$\uC_* (F)$.
Thus, the $l$-th cohomology object of 
$\uC_* (\zeq (Y^\bullet \to W, \BA^n, 0)) (\Spec k)$
equals $\uh^l (\zeq (Y^\bullet \to W, \BA^n, 0)) (\Spec k)$.

\begin{Cor} \label{6CC}
In the situation of Proposition~\ref{6CB}, 
there is a canonical
isomorphism
\[
c_{W / Y}: \uh^l  (\zeq (Y^\bullet \to W, \BA^n, 0)) (\Spec k) \isoto
\Hom (\Mgm (W / Y), \BZ(n)[2n+l]) \; .
\]
Here $\Hom$ denotes morphisms in $\DeffgM$.
The isomorphism is compatible with restriction of $W$ and $Y$
in the sense of Convention~\ref{6conv}.
\end{Cor}

\begin{Proof}
First, we have
\[
\Hom (\Mgm (W / Y), \BZ(n)[2n+l]) =
\Hom_{\DM} (\Mgm (W / Y), \BZ(n)[2n+l])
\]
since the embedding of $\DeffgM$ into $\DM$ is full.
Adjointness of $\otimes$ and $\uHom$ implies that
$\Hom_{\DM} (\Mgm (W / Y), \BZ(n)[2n+l])$ equals
\[
\Hom_{\DM} (\Mgm (\Spec k), \uHom_{\DM} (\Mgm (W / Y), \BZ(n)[2n+l])) \; .
\]
By Proposition~\ref{6CB}, this group  
equals
\[
\Hom_{\DM} (\RC (L (\Spec k)), \RC (\zeq (Y^\bullet \to W, \BA^n, 0))[l]) \; .
\]
Now use the fact that $\RC$ is left adjoint to the inclusion of $\DM$
into $D^- (\ShN)$.
\end{Proof}

For later use, we also note a consequence of the special case $Y = \emptyset$:

\begin{Cor} \label{6Ca}
Let $W \in Sm/k$. Then
the functor on open sub-schemes of $W$
\[
U \longmapsto \uC_* (\zeq (U, \BA^n, 0))
\]
satisfies the Mayer--Vietoris property in the following sense:
given an equality $U = U_1 \cup U_2$ of open sub-schemes of $W$,
the exact sequence 
\[
0 \to 
\zeq (U, \BA^n, 0) \to
\zeq (U_1, \BA^n, 0) \oplus \zeq (U_2, \BA^n, 0) 
\to \zeq (U_1 \cap U_2, \BA^n, 0) 
\]
in $\ShN$ induces an exact triangle
\[ 
\RC (\zeq (U, \BA^n, 0)) \longto
\RC (\zeq (U_1, \BA^n, 0)) \oplus \RC (\zeq (U_2, \BA^n, 0)) \longto 
\]
\[
\longto \RC (\zeq (U_1 \cap U_2, \BA^n, 0)) \longto 
\RC (\zeq (U, \BA^n, 0))[1]
\]
in $\DM$.
\end{Cor}

\begin{Proof}
This follows from
Proposition~\ref{6CB} and 
the Mayer--Vietoris pro\-per\-ty for the functor
$U \mapsto \Mgm (U)$
\cite[Prop.~V.4.1.1]{VSF}.
\end{Proof}

We need to find a way to efficiently generate elements
in the group
\[
\uh^0  (\zeq (Y^\bullet \to W, \BA^n, 0)) (\Spec k) \; .
\] 

\begin{Def} \label{6T}
Let $W \in Sm/k$, and assume that $W$ is of pure dimension $m$.
Fix a closed sub-scheme $Y \subset W$, and an object $V$ in $Sch/k$.
Let $r \ge 0$. Define the sub-sheaf 
\[
\zeq (W \times_k V, m+r)_{Y} \subset \zeq (W \times_k V, m+r)
\]
as follows: $U \in Sm/k$ is mapped to the free Abelian group on the
symbols $(Z)$, where $Z$ runs through those generators of 
$\zeq (W \times_k V, m+r)(U)$ such that for any geometric point
\[
\Spec (\overline{k}) \longto U \times_k Y \longinto U \times_k W 
\]
of $U \times_k Y$, the sub-scheme
\[
Z \times_{U \times_k W} \Spec (\overline{k}) 
\]
of $V_{\overline{k}}$ is empty or of dimension $r$. 
The sheaf of Abelian monoids 
\[
\zeff (W \times_k V, m+r)_{Y} \subset \zeq (W \times_k V, m+r)_{Y}
\]
is defined as the intersection of $\zeq (W \times_k V, m+r)_{Y}$ with
the monoid of effective cycles $\zeff (W \times_k V, m+r)$ in
$\zeq (W \times_k V, m+r)$.
\end{Def}

One checks that $\zeq (W \times_k V, m+r)_{Y}$
and $\zeff (W \times_k V, m+r)_{Y}$
inherit the transfers from $\zeq (W \times_k V, m+r)$. 
If one imposes the defining condition on all geometric points
of $U \times_k W$ instead of just those of $U \times_k Y$,
then one gets $\zeq (W, V, r)(U)$.
Hence $\zeq (W, V, r)$ is a sub-sheaf of 
$\zeq (W \times_k V, m+r)_Y \,$: 
\[
\CD: \zeq (W, V, r) \longinto 
\zeq (W \times_k V, m+r)_{Y} \; .
\] 
Define the natural inclusion 
\[
\iota: \zeq (W \times_k V, m+r)_{Y} \longinto
\zeq (W \times_k V, m+r) \; .
\] 
The Moving Lemma \cite[Thm.~IV.6.3]{VSF} implies
(see \cite[Lemma~IV.6.6]{VSF})
that if both $W$ and $V$ are smooth and projective, then both
$\RC (\iota)$ and $\RC (\iota) \circ \RC (\CD)$
are isomorphisms. Hence $\RC (\CD)$ is an isomorphism if both 
$W$ and $V$ are smooth and projective.
Our aim is to prove this statement under less restrictive
hypotheses on $W$ and $V$. Our result is a variant of
duality for bivariant cycle cohomology
\cite[Thm.~IV.7.4]{VSF}:

\begin{Thm} \label{6dual}
Let $W \in Sm/k$ be quasi-projective, and of pure dimension $m$.
Let $Y \subset W$ be a closed sub-scheme, and $V \in Sch/k$.
Let $r \ge 0$.
Then the inclusion
\[
\CD: \zeq (W, V, r) \longinto 
\zeq (W \times_k V, m+r)_{Y} 
\] 
induces an isomorphism $\RC(\CD)$.
It is compatible with restriction of $W$ and $Y$
in the sense of Convention~\ref{6conv}.
\end{Thm}

\begin{Proof}
We shall follow faithfully the strategy of 
\cite[pp.~172--176]{VSF}. 
Fix compactifications $\oW$ of $W$
and $\oV$ of $V$, with a smooth and projective
$\oW$. For any proper $\oV$-scheme
$\oT$, define the morphism
\[
\alpha_{\oT}: \zeq ({\oW} \times_k {\oT}, m+r) \longto
\zeq ({\oW} \times_k {\oV}, m+r) \longto
\zeq (W \times_k V, m+r)
\]
as the composition of proper push-forward with restriction
\cite[pp.~141--142]{VSF}. Similarly,
define the variant on effective cycle sheaves
\[
\alpha_{\oT}^{eff}: \zeff ({\oW} \times_k {\oT}, m+r) \longto
\zeff (W \times_k V, m+r) \; .
\]
Denote by $\Psi_{\oT} \subset \zeq ({\oW} \times_k {\oT}, m+r)$
the sub-sheaf of Abelian groups generated by 
$(\alpha_{\oT}^{eff})^{-1}(\zeff (W \times_k V, m+r)_{Y})$.
We have
\[
\Psi_{\oV} = \alpha_{\oV}^{-1}(\zeq (W \times_k V, m+r)_{Y}) \; . 
\]
We claim that the inclusion  
\[
\iota_{\oT}: \Psi_{\oT} \longinto \zeq ({\oW} \times_k {\oT}, m+r) \; .
\]
induces a quasi-isomorphism $\uC_* (\iota_{\oT})$. In order to prove this
claim, imitate the proof of \cite[Prop.~IV.7.3]{VSF}.
One uses the Moving Lemma we already cited.
It is here that the projectivity assumption on $\oW$ enters.

Next, one imitates the proof of \cite[Thm.~IV.7.4]{VSF}, using 
the above instead of \cite[Prop.~IV.7.3]{VSF}, to see that the inclusion
\[
\iota : \zeq (W \times_k V, m+r)_{Y}  
\longinto \zeq (W \times_k V, m+r) 
\]
induces a quasi-isomorphism $\uC_* (\iota)$.
The same observation, applied to the case where $Y=W$ implies that
the composition $\RC(\iota) \circ \RC(\CD)$ is an isomorphism.
(This is of course the original statement of \cite[Thm.~IV.7.4]{VSF}).
\end{Proof}

We unite the assumptions from \ref{6CB} and \ref{6T}:
$W \in Sm/k$ is of pure dimension $m$, and
$Y \subset W$ is a closed sub-scheme
such that arbitrary intersections of
the components of $Y$ are smooth. 
Observe that the condition on elements in  
$\zeq (W \times_k \BA^n, m)_Y$ ensures in particular
that they intersect properly with
the components $Y_j \times_k \BA^n$ of $Y \times_k \BA^n$. 
This allows to define an inverse image $\delta$
making the following diagram commutative:
\[
\xymatrix@R-20pt{
\zeq (W, \BA^n, 0) \ar[dd]_{\delta} \ar@{^{ (}->}[r]^-{\CD} & 
\zeq (W \times_k \BA^n, m)_Y \ar[dd]^{\delta} \\
& \\
\bigoplus_j \zeq (Y_j, \BA^n, 0) \ar@{^{ (}->}[r]^-{\CD}  & 
\bigoplus_j \zeq (Y_j \times_k \BA^n, \dim Y_j)
\\}
\]
The group in the lower right corner has to be modified if
$Y_j$ has several components of different dimension. 
But in fact, $\zeq (W \times_k \BA^n, m)_Y$ is defined such that
$\delta$ not only exists but maps
$\zeq (W \times_k \BA^n, m)_Y$ to 
$\oplus_j \zeq (Y_j, \BA^n, 0)$.
Thus we may enlarge the complex 
\[
\zeq (Y^\bullet \to W,\BA^n,0) 
\]
by replacing its zeroth component $\zeq (W, \BA^n, 0)$
with $\zeq (W \times_k \BA^n, m)_Y\, $:
\[
\zeq (Y^\bullet \to W,\BA^n,0) \subset
\zeq (Y^\bullet \to W,\BA^n,0)' \; .
\]
When $W$ is quasi-projective, then Theorem~\ref{6dual}
shows that this inclusion of complexes
induces an isomorphism after application of $\RC$. In the ge\-ne\-ral case,
there is a (mainly notational) complication since we do not know
whether the functor on open sub-schemes of $W$
\[
U \longmapsto \uC_* (\zeq (U \times_k \BA^n, m)_Y)
\]
satisfies the Mayer--Vietoris property. Therefore, we fix an
additional geometric datum, namely a finite open covering $\{ W_\alpha \}$
of $W$ by quasi-projective schemes. 
Consider the covering of $Y$ induced by
$\{ W_\alpha \}$, and define 
\[
\zeq (Y^\bullet \to W,\BA^n,0)'' 
\]
as the simple complex of the double \v{C}ech complex
associated to this covering and $\zeq (Y^\bullet \to W,\BA^n,0)'$. 
For example, the components
of degree $0$ and $1$ are
\[
\bigoplus_\alpha \zeq (W_\alpha \times_k \BA^n, m)_Y 
\]
and
\[
\bigoplus_{\alpha,j} \zeq (W_\alpha \cap Y_j, \BA^n, 0) \oplus
\bigoplus_{\alpha_1 \ne \alpha_2} 
        \zeq ((W_{\alpha_1} \cap W_{\alpha_2}) \times_k \BA^n, m)_Y \; .
\]
Consider the natural inclusion
\[
\CD: \zeq (Y^\bullet \to W, \BA^n, 0) \longinto
\zeq (Y^\bullet \to W, \BA^n, 0)' \; ,
\]
and use the same symbol for the composition of $\CD$
with the co-augmentation
\[
\zeq (Y^\bullet \to W, \BA^n, 0)' \longto
\zeq (Y^\bullet \to W, \BA^n, 0)'' \; .
\]

\begin{Cor} \label{6Ia}
The morphism
\[
\CD : \zeq (Y^\bullet \to W, \BA^n, 0) \longto 
\zeq (Y^\bullet \to W, \BA^n, 0)'' 
\] 
induces an isomorphism $\RC (\CD)$.
In particular, different choices of coverings $W = \cup_\alpha W_\alpha$
give rise to the same object $\RC(\zeq (Y^\bullet \to W, \BA^n, 0)'')$.
The isomorphism $\RC (\CD)$ is compatible with restriction of $W$ and $Y$
in the sense of Convention~\ref{6conv}.
\end{Cor}

\begin{Proof}
This follows from Theorem~\ref{6dual} and Corollary~\ref{6Ca}.
\end{Proof}

\begin{Def} \label{6Ib}
Define the sub-sheaf with transfers
\[
\zeq (W, m-n)_Y \subset \zeq (W, m-n)
\]
as the sub-sheaf of cycles having empty intersection with $Y$.
\end{Def}

Note that $\zeq (W, m-n)_Y$
behaves contravariantly with respect to 
restriction of $W$ and $Y$.
Flat pull-back defines a morphism 
\[
p_{\BA^n}^*: \zeq (W, m-n)_Y [0] \longto \zeq (Y^\bullet \to W, \BA^n, 0)' \; .
\]
Its composition with the co-augmentation
\[
\zeq (Y^\bullet \to W, \BA^n, 0)' \longto
\zeq (Y^\bullet \to W, \BA^n, 0)'' 
\]
induces a morphism on the level of $\uh^0 (\Spec k)$,
denoted by the same symbol $p_{\BA^n}^*$
Putting everything together, we obtain:

\begin{Thm} \label{6L}
Let $W \in Sm/k$ be of pure dimension $m$, and
$Y \subset W$ a closed sub-scheme
such that arbitrary intersections of
the components of $Y$ are smooth.  
Then there is a unique morphism
\[
cyc_{W / Y}: \uh^0 (\zeq (W, m-n)_Y) (\Spec k) \longto  
\Hom (\Mgm (W / Y), \BZ(n)[2n]) 
\]
making the following diagram commute:
\[
\xymatrix@R-20pt{
\uh^0 (\zeq (W, m-n)_Y) (\Spec k) \ar[r]^-{p_{\BA^n}^*} \ar[dd]_{cyc_{W / Y}} &
\uh^0 (\zeq (Y^\bullet \to W, \BA^n, 0)'') (\Spec k) \\
& \\
\Hom (\Mgm (W / Y), \BZ(n)[2n]) &
\uh^0 (\zeq (Y^\bullet \to W, \BA^n, 0)) (\Spec k) 
          \ar[l]_-{c_{W / Y}}^-{\cong} \ar[uu]_{\CD}^{\cong}
\\}
\]
Here $\Hom$ denotes morphisms in $\DeffgM$.
The morphism $cyc_{W / Y}$ is compatible with restriction of $W$ 
and $Y$ in the sense of Convention~\ref{6conv}.
\end{Thm}

\begin{Proof}
Apply Corollaries~\ref{6CC} and \ref{6Ia}. 
\end{Proof}

\begin{Rem} \label{6M}
Another type of compatibility property with respect to change of
$W$ and $Y$ is useful. Assume given a second pair $Y_1 \subset W_1$
of schemes satisfying the hypotheses of Theorem~\ref{6L}, and a
proper morphism $W_1 \to W$ identifying $Y_1$ with the fibre product
$W_1 \times_W Y$, and inducing an isomorphism from $W_1 - Y_1$ to
$W - Y$. Theorem~\ref{3aB} tells us that 
\[
\Mgm (W_1 / Y_1) \longto \Mgm (W / Y)
\]
is an isomorphism. On the other hand, we clearly have
\[
\zeq (W, m-n)_Y = \zeq (W_1, m-n)_{Y_1} \; .
\]
It is easy to see that in this situation, the diagram
\[
\xymatrix@R-20pt{
\uh^0 (\zeq (W, m-n)_Y) (\Spec k) \ar[r]^-{=} \ar[dd]_{cyc_{W / Y}} &
\uh^0 (\zeq (W_1, m-n)_{Y_1}) (\Spec k) \ar[dd]_{cyc_{W_1 / Y_1}} \\
& \\
\Hom (\Mgm (W / Y), \BZ(n)[2n]) \ar[r]^-{\cong} &
\Hom (\Mgm (W_1 / Y_1), \BZ(n)[2n]) 
\\}
\]
commutes.
\end{Rem}

When $Y$ is empty, then \cite[Cor.~V.4.2.5]{VSF}
tells us that $cyc_W := cyc_{W / \emptyset}$ is an isomorphism. 
We have not tried to see whether the analogous statement for 
$cyc_{W / Y}$ is true when $Y$ is non-empty. \\

\begin{Proofof}{Theorem~\ref{6main}}
We start by introducing the notation
$Cyc^n (W)_Y$ for the group $\zeq (W, m-n)_Y (\Spec k)$, when $W$ and
$Y \subset W$ are as before. By definition,
$Cyc^n (W)_Y$ is the group of codimension $n$ cycles
on $W$ not meeting $Y$. Write $Cyc^n (W)$ for the group of all
codimension $n$ cycles on $W$.
Our proof relies on the following principles, which are
consequences of Theorem~\ref{6L},
and of the definition of the tensor structure on $\DM$:
\begin{itemize}
\item[(A)] For $i = 1, 2$, let $W_i \in Sm/k$ be of pure dimension, with 
closed sub-schemes $Y_i \subset W_i$,   
such that arbitrary intersections of
the components of the $Y_i$ are smooth.
Then any 
\[ 
c \in Cyc^n (W_1 \times_k W_2)_{Y_1 \times_k W_2 \cup W_1 \times_k Y_2}
\]
defines a pairing
\[
(\argdot,\argdot)_c: \Mgm (W_1 / Y_1) \otimes \Mgm (W_2 / Y_2)  
                                       \longto \BZ (n)[2n] \; ,
\]
or equivalently, a morphism
\[
\varepsilon_c: \Mgm (W_1 / Y_1) \longto \Mgm (W_2 / Y_2)^*(n)[2n] 
\]
in $\DgM$. The morphism $\varepsilon_c$ is induced by a morphism
of Nisnevich sheaves 
\[
e_{c'}: c(\argdot, W_1) /  c(\argdot, Y_1) \longto 
\zeq (Y_2^\bullet \to W_2, \BA^n, 0)
\]
defined as follows: use Corollary~\ref{6Ia} to move the pull-back
$p_{\BA^n}^* (c)$ into a cycle 
\[
c' \in \zeq (W_1 \times_k W_2, \BA^n, 0)(\Spec k) 
\]
having empty intersection with
$Y_1 \times_k W_2 \cup W_1 \times_k Y_2$.
Given $U \in Sm/k$
and $\CZ \in c(U, W_1)$, pull back $\CZ$ to $W_2$, giving
$\CZ_{W_2} \in c(U \times_k W_2, W_1)$. Similarly,
pull back $c'$ to $U$, giving $c'_U$. Now consider the cycle
\[
c' \cap \CZ_{W_2} := Cor_{W_1 \times_k W_2 \times_k \BA^n / W_1 \times_k W_2}
(c'_U \otimes \CZ_{W_2}) \; ,
\]
where $Cor_{W_1 \times_k W_2 \times_k \BA^n / W_1 \times_k W_2}$
is the correspondence homomorphism from \cite[II.3.7]{VSF}.
By \cite[Cor.~II.3.7.5]{VSF}, we have
\[
c' \cap \CZ_{W_2} \in \zeq (W_2, W_1 \times_k \BA^n, 0)(U) \; .
\]
Furthermore, finiteness of $\CZ$ over $U$ implies that $c' \cap \CZ_{W_2}$ is
finite over $U \times_k W_2 \times_k \BA^n$. Push-forward via the
projection $p_2$ to this product then yields
\[
e_{c'} (\CZ) := p_{2 \, *} (c' \cap \CZ_{W_2}) \in \zeq (W_2, \BA^n, 0)(U) \; .
\]
In fact, 
$e_{c'} (\CZ)$ lies in $\zeq (Y_2^\bullet \to W_2, \BA^n, 0)(U)$
since $c'$ has empty intersection with
$W_1 \times_k Y_2$. Furthermore, $e_{c'} (\CZ)$ only depends on
the class of $\CZ$ modulo $c(U, Y_1)$ since $c'$ has
empty intersection with $Y_1 \times_k W_2$.
\item[(B)] Let $W \in Sm/k$ be of pure dimension, with a
closed sub-scheme $Y \subset W$,   
such that arbitrary intersections of
the components of $Y$ are smooth. Let  
\[ 
c \in Cyc^n (W \times_k W)_{Y \times_k W \cup W \times_k Y} \; .
\]
Then a sufficient condition for the pairing 
\[
(\argdot,\argdot)_c: \Mgm (W / Y) \otimes \Mgm (W / Y)  
                                       \longto \BZ (n)[2n] 
\]
to be symmetric,
or equivalently, for the morphism
\[
\varepsilon_c: \Mgm (W / Y) \longto \Mgm (W / Y)^*(n)[2n] 
\]
to be auto-dual: $\varepsilon_c = \varepsilon_c^*(n)[2n]$,
is the symmetry of the cycle $c$ 
in $Cyc^n (W \times_k W)$.
\item[(C)] For $i = 1, \ldots, 4$, let 
$W_i \in Sm/k$ be of pure dimension, with 
closed sub-schemes $Y_i \subset W_i$,   
such that arbitrary intersections of
the components of the $Y_i$ are smooth. For $i = 1, 2$, let 
$j_i: W_i \into W_{i+2}$ be an open immersion mapping $Y_i$ into
$Y_{i+2}$. Let 
\[ 
c_{1,4} \in Cyc^n (W_1 \times_k W_4)_{Y_1 \times_k W_4 \cup W_1 \times_k Y_4}
\]
and
\[ 
c_{3,2} \in Cyc^n (W_3 \times_k W_2)_{Y_3 \times_k W_2 \cup W_3 \times_k Y_2} 
\; .
\]
Then a sufficient condition for the diagram
\[
\xymatrix@R-20pt{
\Mgm (W_1 / Y_1) \ar[r]^-{\varepsilon_{c_{1,4}}} \ar[dd]_{j_1} &
\Mgm (W_4 / Y_4)^*(n)[2n] \ar[dd]^{j_2^*(n)[2n]} \\
& \\
\Mgm (W_3 / Y_3) \ar[r]^-{\varepsilon_{c_{3,2}}} &
\Mgm (W_2 / Y_2)^*(n)[2n]  
\\}
\]
to commute is the equality of cycles
\[
(j_1,\id_{W_2})^* c_{3,2} = (\id_{W_1},j_2)^* c_{1,4}
\]
in $Cyc^n (W_1 \times_k W_2)$.
\end{itemize} 
Furthermore, using the compatibility of moving cycles with 
correspondence homomorphisms and direct images, and 
\cite[Prop.~II.3.7.6]{VSF}, one sees:
\begin{itemize}
\item[(D)] For $i = 1, 2$, let $W_i \in Sm/k$ be of pure dimension $m_i$, with a
closed sub-scheme $Y_1 \subset W_1$,   
such that arbitrary intersections of
the components of $Y_1$ are smooth. Assume that $m_2 \ge n$, and that
\[ 
c \in Cyc^n (W_1 \times_k W_2)_{Y_1 \times_k W_2} \cap
\zeq (W_1, W_2, m_2 - n)(\Spec k) \; .
\]
Then the morphism
\[
\varepsilon_c: \Mgm (W_1 / Y_1) \longto \Mgm (W_2)^*(n)[2n] 
\]
associated to $c$
is also induced by the composition of the morphism of Nisnevich sheaves 
\[
e'_c: c(\argdot, W_1) /  c(\argdot, Y_1) \longto 
\zeq (W_2 \times_k \BA^n, m_2) 
\]
with the inverse of the isomorphism $\RC(\CD)$ from Theorem~\ref{6dual},
where $e'_c$ is defined as follows: 
given $U \in Sm/k$
and $\CZ \in c(U, W_1)$, pull back $\CZ$ to $\BA^n$, giving
$\CZ_{\BA^n} \in c(U \times_k \BA^n, W_1)$. Similarly,
pull back $c$ to $U \times_k \BA^n$, giving $c_{U \times_k \BA^n}$. 
Now consider the cycle
\[
c \cap \CZ_{\BA^n} := Cor_{W_1 \times_k W_2 \times_k \BA^n / W_1 \times_k \BA^n}
(c_{U \times_k \BA^n} \otimes \CZ_{\BA^n}) \; .
\]
By \cite[Cor.~II.3.7.5]{VSF}, we have
\[
c \cap \CZ_{\BA^n} \in \zeq (\BA^n, W_1 \times_k W_2, m_2-n)(U) \; .
\]
Furthermore, finiteness of $\CZ$ over $U$ implies that $c \cap \CZ_{\BA^n}$ is
finite over $U \times_k W_2 \times_k \BA^n$. Push-forward via the
projection $p_2$ to this product then yields
\[
e'_c (\CZ) := p_{2 \, *} (c \cap \CZ_{\BA^n}) 
\in \zeq (\BA^n, W_2, m_2-n)(U) \; ,
\]
and the latter group is contained in
$\zeq (W_2 \times_k \BA^n, m_2)(U)$.
Observe that $e'_c (\CZ)$ only depends on
the class of $\CZ$ modulo $c(U, Y_1)$ since $c$ has 
empty intersection with $Y_1 \times_k W_2$.
\end{itemize}
Choose and fix a smooth compactification $\oX$ of $X$ such that
$\doX := \oX - X$ is a normal crossing divisor with smooth
irreducible components. By \cite[Prop.~V.4.1.5]{VSF}, we have
a canonical isomorphism between $\Mgm (\oX / \doX)$ and $\Mcgm (X)$. 
Applying principle~(A), we see that the diagonal
\[
\Delta \in Cyc^n ({\oX} \times_k X)_{{\doX} \times_k X}
\]
induces a morphism
\[
\varepsilon_\Delta: \Mcgm (X) \longto \Mgm (X)^*(n)[2n] \; .
\] 
Principle~(D) and the proof of \cite[Thm.~V.4.3.7]{VSF} show that
this is the morphism $\mu_X$. To say that 
\[
\xymatrix@R-20pt{
\Mgm (X) \ar[rr]^-{\nu_X = \mu_X^*(n)[2n]} \ar[dd]_{\beta} &&
\Mcgm (X)^*(n)[2n] \ar[dd]^{\beta^*(n)[2n]} \\
&& \\
\Mcgm (X) \ar[rr]^-{\mu_X} &&
\Mgm (X)^*(n)[2n]  
\\}
\]
commutes is equivalent to saying that
the pairing 
\[
\Mgm (X) \otimes \Mgm (X) \stackrel{\id \otimes \beta} \longto
\Mgm (X) \otimes \Mcgm (X) \stackrel{(\argdot,\argdot)_\Delta} \longto
\BZ (n) [2n]
\]
is symmetric. By principle~(B), this is indeed the case since the
restriction of the cycle $\Delta$ to $X \times_k X$ is symmetric.

Now for the construction of 
\[
\eta_X: \dMgm (X) \longto \dMgm (X)^*(n)[2n-1] \; .
\]  
Recall that by Proposition~\ref{2C}, there is a canonical isomorphism
\[
\cone \left( \Mgm \left( X \coprod {\doX} \right) 
\to \Mgm (\oX) \right) \isoto \dMgm (X) \; . 
\]
From this, one deduces that $\dMgm (X) \otimes \dMgm (X) [1]$ is represented
by the complex $\uC_* (s \CL_*)$, where $\CL_*$ is the complex
concentrated in degrees $-1$, $0$, and $1$
\[
\xymatrix@R-20pt{
L(X \times_k X) \oplus L(X \times_k {\doX}) \oplus 
L({\doX} \times_k X) \oplus L({\doX} \times_k {\doX}) \ar[dd] \\
& \\
L({\oX} \times_k X) \oplus L(X \times_k {\oX}) \oplus 
L({\oX} \times_k {\doX}) \oplus L({\doX} \times_k {\oX}) \ar[dd] \\
& \\
L({\oX} \times_k {\oX})
\\}
\]
the differentials being induced by the inclusions. In particular,
one sees that the complex $\CL'_*$ 
concentrated in degrees $-1$ and $0$
\[
\xymatrix@R-20pt{
L(X \times_k X) \oplus L(X \times_k {\doX}) \oplus 
L({\doX} \times_k X) \ar[dd] \\
& \\
L({\oX} \times_k X) \oplus L(X \times_k {\oX}) 
\\}
\]
is a quotient of $\CL_*$. By the Mayer--Vietoris property for 
the functor $\Mgm$ \cite[Prop.~V.4.1.1]{VSF}, $\uC_*(\CL'_*)$
is canonically quasi-isomorphic to the complex $\uC_*(\CL''_*)$,
where $\CL''_*$ is the complex  
concentrated in degrees $-1$ and $0$
\[
\xymatrix@R-20pt{
L({\oX} \times_k {\oX} - {\doX} \times_k {\doX} - X \times_k X)  \ar[dd] \\
& \\
L({\oX} \times_k {\oX} - {\doX} \times_k {\doX}) 
\\}
\]
This shows that there is a canonical morphism from 
$\dMgm (X) \otimes \dMgm (X) [1]$ to
\[
\Mgm (({\oX} \times_k {\oX} - {\doX} \times_k {\doX}) / 
({\oX} \times_k {\oX} - {\doX} \times_k {\doX} - X \times_k X)) \; .
\]
Applying principle~(A), we see that the diagonal
\[
\Delta \in 
Cyc^n ({\oX} \times_k {\oX} - {\doX} \times_k {\doX})_{
{\oX} \times_k {\oX} - {\doX} \times_k {\doX} - X \times_k X}
\]
induces a morphism
\[
\varepsilon_\Delta: \dMgm (X) \longto \dMgm (X)^*(n)[2n-1] \; .
\] 
We define this to be the morphism $\eta_X$. Principle~(B) shows that
$\eta_X$ is auto-dual. In order to see that $\eta_X$ 
fits into a morphism of exact triangles
\[
\xymatrix@R-20pt{
\dMgm (X) \ar[r]^-{\eta_X} \ar[dd]_{\alpha} &
\dMgm (X)^*(n)[2n-1] \ar[dd]^{\gamma^*(n)[2n]} \\
& \\
\Mgm (X) \ar[r]^-{\nu_X} \ar[dd]_{\beta} &
\Mcgm (X)^*(n)[2n] \ar[dd]^{\beta^*(n)[2n]} \\
& \\
\Mcgm (X) \ar[r]^-{\mu_X} \ar[dd]_{\gamma} &
\Mgm (X)^*(n)[2n] \ar[dd]^{\alpha^*(n)[2n]} \\
& \\
\dMgm (X)[1] \ar[r]^-{\eta_X[1]} &
\dMgm (X)^*(n)[2n]
\\}
\]
apply principle~(C). By \cite[Thm.~V.4.3.7~3.]{VSF},
$\mu_X$ and $\nu_X$ are isomorphisms, and hence so is $\eta_X$.
In order to check that $\eta_X$ does not depend on the choice of
the smooth compactification $\oX$, use Remark~\ref{6M}
together with the fact that the system
of such compactifications is filtering.
\end{Proofof}

Fix a compactification $\oX$ of $X$ (which may be non-smooth), 
and set $\doX := \oX - X$.
Write $j$ for the open immersion of $X$, 
and $i_{\doX}$ for the closed immersion of $\doX$ into $\oX^m$.  
We are thus in the situation considered in Section~\ref{3},
with $Y = Y' = \doX$ and $W = \oX$.
By Proposition~\ref{2C}, we have a canonical isomorphism
\[
\Mgm (\doX, \iXjZ) \isoto \dMgm (X)[1] \; .
\] 
By duality, this gives
\[
\alpha_X: \dMgm (X)^*[-1] \isoto \Mgm (\doX, \iXjZ)^* \; .
\]
Preceding $\alpha_X(n)[2n]$ with
the auto-duality $\eta_X$ from Theorem~\ref{6main},
we obtain:

\begin{Cor} \label{6N}
There is a canonical isomorphism
\[
\dMgm (X) \isoto \Mgm (\doX, \iXjZ)^*(n)[2n] \; .
\]
\end{Cor}


\bigskip

%
%

\section{Localization}
\label{7}



Throughout this section, we assume that $k$ admits resolution of singularities.
Fix closed immersions $Y \into Y' \into W$ in $Sch / k$. Write
$j$ for the open immersion of $W-Y'$, and $i_Y$ for the closed
immersion of $Y$ into $W$.

\begin{Def} \label{7A}
Assume that $W - Y' \in Sm/k$.  \\[0.1cm]
(a)~If $W - Y'$ is of pure dimension $n$, we put
\[
\Mcgm (Y, \iastYjZ) := \Mgm(Y, \iYjZ)^*(n)[2n] \; .
\]
(b)~If $W - Y' = \coprod_\alpha W_\alpha$ is the decomposition
of $W - Y'$ into connected components, then the 
\emph{motive with compact support of $Y$ and with
coefficients in $\iastYjZ$} is defined as
\[
\Mcgm (Y, \iastYjZ) := \bigoplus_\alpha \Mcgm (Y, \iastYjaZ) \; .
\]
Here, the $j_\alpha$ denote the open immersions of the $W_\alpha$ into $W$.
\end{Def}

Using \cite[Cor.~V.4.3.6]{VSF} and the definition of $\Mgm(Y, \iYjZ)$ 
(Definition~\ref{3A}), one sees that 
$\Mcgm (Y, \iastYjZ)$ belongs to $\DeffgM$. 

\begin{Rem} \label{7B}
The object on the right hand side in Definition~\ref{7A} is defined
without the hypothesis of smoothness on $W - Y'$. In general, the object
$\Mgm(Y, \iYjZ)$ is dual to what should be considered as
the motive with compact support of $Y$ and with
coefficients in $\iastYjDZ$, where $\BD (\BZ)$ is the coefficient system
on $W - Y'$ which is Poincar\'e-dual to $\BZ$. 
\end{Rem}

Now assume given a filtration
\[
\emptyset = \FF_{-1} Y \subset \FF_0 Y \subset \ldots \subset 
\FF_d Y = Y
\]
of $Y$ by closed sub-schemes. It induces a stratification of $Y$ by
locally closed sub-schemes $Y_m := \FF_m Y - \FF_{m-1} Y$, for
$m = 0, \ldots, d$. Define $W^m$ as the complement of $\FF_{m-1} Y$ in
$W$. 
Write $i_{Y_m}$ for the closed immersion of $Y_m$ into $W^m$.
By abuse of notation, we use the letter $j$ to denote also the open
immersions of $W - Y'$ into $W^m$. 

\begin{Thm}[Localization] \label{7main}
Assume that $W - Y' \in Sm/k$.
Then there is a cano\-ni\-cal chain of morphism
\[
M_0 = \Mcgm (Y, \iastYjZ)
\stackrel{\gamma_0}{\longto} M_1 
\stackrel{\gamma_1}{\longto} M_2  
\stackrel{\gamma_2}{\longto} \ldots
\stackrel{\gamma_{d-1}}{\longto} M_d
\stackrel{\gamma_d}{\longto} M_{d+1} = 0
\]
in $\DgM$. For each 
$m \in \{ 0, \ldots, d\}$, there is a canonical exact triangle
\[
\Mcgm (Y_m, \iastYmjZ) \longto M_m 
\stackrel{\gamma_m}{\longto} M_{m+1} \longto \Mcgm (Y_m, \iastYmjZ)[1] 
\]
in $\DgM$. In particular, all the $M_m$ are in $\DeffgM$.
\end{Thm}
       
\begin{Proof}
Dualize co-localization (Theorem~\ref{3C}).
\end{Proof}

\begin{Cor} \label{7C}
In the above situation, assume that $Y = \doX := \oX - X$,
with $X \in Sm/k$, and $\oX$ 
a compactification of $X$ (which may be non-smooth).
Write $(\doX)_m:= Y_m$.  
Then there is a cano\-ni\-cal chain of morphisms
\[
M_0 = \dMgm (X)
\stackrel{\gamma_0}{\longto} M_1 
\stackrel{\gamma_1}{\longto} M_2  
\stackrel{\gamma_2}{\longto} \ldots
\stackrel{\gamma_{d-1}}{\longto} M_d
\stackrel{\gamma_d}{\longto} M_{d+1} = 0
\]
in $\DgM$. For each 
$m \in \{ 0, \ldots, d\}$, there is a canonical exact triangle
\[
\Mcgm ((\doX)_m, \iastXmjZ) \longto M_m 
\stackrel{\gamma_m}{\longto} M_{m+1} \longto \Mcgm ((\doX)_m, \iastXmjZ)[1] 
\]
in $\DgM$. In particular, all the $M_m$ are in $\DeffgM$.
\end{Cor}

\begin{Proof}
This follows from Corollary~\ref{6N} and
Theorem~\ref{7main}.
\end{Proof}

\begin{Rem} \label{7D}
Given our definition, it is easy to deduce from duality and
Theorems~\ref{3aA}
and \ref{4A} that $\Mcgm (Y, \iastYjZ)$ is invariant under abstract
blow-up and under analytical isomorphism. In concrete situations,
this observation helps to control the ``graded pieces''
$\Mcgm ((\doX)_m, \iastXmjZ)$ of $\dMgm (X)$.
\end{Rem}

\bigskip

%
%

\section{The case of normal crossings}
\label{8}



Throughout this section, we consider the following situation: 
$X$ lies in $Sm/k$ and is irreducible, 
and $\oX \in Sm/k$ is a 
smooth compactification of $X$, such that $\doX := \oX - X$ is a normal
crossing divisor with smooth irreducible components. 
We stratify $\doX$ by defining
$(\doX)^m$ as the geometric locus of
points lying on exactly $m$ irreducible components, $m = 1, \ldots, \dim X$.
Note that the $(\doX)^m$ are all smooth.
Denote by $j$ the open immersion of $X$, and by $i_m$ the immersion of
$(\doX)^m$ into $\oX$.
We re-state co-localization and localization (Corollaries~\ref{3D} and \ref{7C}):

\begin{Cor} \label{8A}
Denote by $n$ the dimension of $X$,
and assume that $k$ admits resolution of singularities. \\[0.1cm]
(a)~There is a cano\-ni\-cal chain of morphisms
\[
M^n = 0 \stackrel{\gamma^{n-1}}{\longto} M^{n-1}  
\stackrel{\gamma^{n-2}}{\longto} M^{n-2}
\stackrel{\gamma^{n-3}}{\longto} \ldots
\stackrel{\gamma^0}{\longto} M^0 = \dMgm (X)
\]
in $\DeffgM$. 
For each 
$m \in \{ 1, \ldots, n\}$, there is a canonical exact triangle
\[
\Mgm ((\doX)^m, \iXnmjZ) [-1] \longto M^{n-m+1} 
\stackrel{\gamma^{n-m}}{\longto} M^{n-m} \longto \Mgm ((\doX)^m, \iXnmjZ)
\]
in $\DeffgM$. \\[0.1cm]
(b)~There is a cano\-ni\-cal chain of morphisms
\[
M_0 = \dMgm (X)
\stackrel{\gamma_0}{\longto} M_1 
\stackrel{\gamma_1}{\longto} M_2  
\stackrel{\gamma_2}{\longto} \ldots
\stackrel{\gamma_{n-2}}{\longto} M_{n-1}
\stackrel{\gamma_{n-1}}{\longto} M_n = 0
\]
in $\DeffgM$. For each 
$m \in \{ 1, \ldots, n\}$, there is a canonical exact triangle
\[
\Mcgm ((\doX)^m, \iastXnmjZ) \longto M_{n-m} 
\stackrel{\gamma_{n-m}}{\longto} M_{n-m+1} \longto \Mcgm ((\doX)^m, \iastXnmjZ)[1] 
\]
in $\DeffgM$.
\end{Cor}

The aim of this section is give a description of the motives
with coefficients $\Mgm ((\doX)^m, \iXnmjZ)$ and $\Mcgm ((\doX)^m, \iastXnmjZ)$
occurring in the above statement.
Denote by $N_m$ the normal bundle of $(\doX)^m$ in $\oX$. For any component 
$D$ of the divisor $\doX$, consider 
the normal bundle of $(\doX)^m \cap D$ in $D$. 
Using all possible $D$,
this gives $m$ sub-bundles of $N_m$ of codimension one. Their intersection is the zero bundle over $(\doX)^m$.
Define $N_m^*$ as the complement in $N_m$ of 
the union of these sub-bundles. 
Fix an order $\prec$ of the index set of components of $\doX$, i.e., write
\[
\doX = \bigcup_{i=1}^r D_i \; ,
\]
where $D_1, \ldots, D_r$ are the $r$ distinct components of $\doX$.
Note that $N_m^*$ is a torsor under a torus of dimension $m$, and that
the choice of $\prec$ gives an identification of this torus with $\BG_m^m$. 

\begin{Thm} \label{8main}
(a)~The order $\prec$ induces an isomorphism 
\[
\Mgm ((\doX)^m, \iXnmjZ) \isoto \Mgm (N_m^*) [m]
\]
in $\DeffgM$. \\[0.1cm]
(b)~The order $\prec$ induces an isomorphism 
\[
\Mcgm ((\doX)^m, \iastXnmjZ) \isoto \Mcgm (N_m^*) [-m]
\]
in $\DeffgM$. 
\end{Thm}

Note that we do not assume resolution of singularities for Theorem~\ref{8main}.
Our main computational tool in its proof will be the following:

\begin{Lemma} \label{8B}
Let $S \in Sch/k$, and $V$ a vector bundle of rank $m$ over $S$.
Assume given sub-bundles $V_i$, $i = 1, \ldots, \ell$ 
of codimension one of $V$, with $1 \le \ell \le m$.
Assume also that Zariski-locally over $S$, there exist trivializations 
$V \cong \BA_S^m$ identifying $V_i$ with the hyperplane given by the
vanishing of the $i$-th coordinate, $i = 1, \ldots, \ell$.
Define $V^*$ as the complement in $V$ of the union of the $V_i$. 
Then for any proper subset $I$ of $\{1, \ldots, \ell \}$, the immersion
$V - V^* - \cup_{i \in I} V_i \into V - \cup_{i \in I} V_i$ induces an isomorphism
\[
\Mgm (V - V^* - \cup_{i \in I} V_i) \isoto \Mgm (V - \cup_{i \in I} V_i)
\]
in $\DeffgM$. 
\end{Lemma}

\begin{Proof}
The Mayer--Vietoris property
\cite[Prop.~V.3.1.3]{VSF} for the functor $L$
shows that we may assume that
a trivialization as in the hypotheses exists globally.
Then the immersion $V^* \into V - \cup_{i \in I} V_i$
is isomorphic to 
\[
\BG_{m,S}^\ell \times_S \BA_S^{m-\ell} \longinto 
\BG_{m,S}^p \times_S \BA_S^{m-p} \; ,
\]
for some $p < \ell$.
Hence $V - V^* - \cup_{i \in I} V_i \into V - \cup_{i \in I} V_i$ is a homotopy equivalence 
(in fact, both sides are homotopic to $\BG_{m,S}^p$).
\end{Proof}

\begin{Proofof}{Theorem~\ref{8main}}
(a)~Recall that $\Mgm ((\doX)^m, \iXnmjZ)$ is defined as the motive $\Mgm (\FY)$, where
$\FY$ is the diagram
\[
\xymatrix@R-20pt{
\doX' - (\doX)^m \ar[r] \ar[dd] & 
\oX' - (\doX)^m \ar[dd] \\
& \\
\doX' \ar[r] & 
\oX' 
\\}
\]
and $\oX'$ and $\doX'$ denote the complements of the strata $(\doX)^p$, for $p > m$.
Now because of our assumptions on $\oX$ and $\doX$, this diagram is analytically
isomorphic to the following diagram $\FN$, in the sense that the assumptions of Theorem~\ref{4A}
are satisfied:
\[
\xymatrix@R-20pt{
N_m - N_m^* - 0(\doX)^m \ar[r] \ar[dd] & 
N_m - 0(\doX)^m \ar[dd] \\
& \\
N_m - N_m^* \ar[r] & 
N_m 
\\}
\]
Here, $0(\doX)^m$ denotes the zero section of the bundle $N_m$ over $(\doX)^m$.
Analytical invariance thus allows us to compute $\Mgm (\FN)$
instead of $\Mgm (\FY)$.
In order to do so, 
consider the upper line 
\[
N_m - N_m^* - 0(\doX)^m \longto
N_m - 0(\doX)^m
\]
of $\FN$, and the finite covering 
\[
N_m - 0(\doX)^m = \bigcup_{i=1}^m N_m - V_i \; .
\]
Here, the $V_i$ are the sub-bundles of codimension one constructed from
the components of $\doX$.
For $1 \le i_1 < \ldots < i_p \le m$,
we have:
\[
\bigcap_{q=1}^p N_m - V_{i_q} = N_m - \bigcup_{q=1}^p V_{i_q} \; . 
\]
Define the double complex $L(\FN')$ of Nisnevich sheaves 
as follows: first consider the \v{C}ech complex associated to the
above covering, and to the induced covering of 
$N_m - N_m^* - 0(\doX)^m$, then add the lower line of $L(\FN)$: 
\[
\xymatrix@R-20pt{
0 = L(\emptyset) \ar[r] \ar[dd] &
L(N_m^*) \ar[dd] \\
& \\
\vdots \ar[dd] &
\vdots \ar[dd] \\
& \\
\prod_{i_1 < i_2} 
L(N_m - N_m^* - (V_{i_1} \cup V_{i_2})) \ar[r] \ar[dd] &
\prod_{i_1 < i_2} 
L(N_m - (V_{i_1} \cup V_{i_2})) \ar[dd] \\
& \\
\prod_i L(N_m - N_m^* - V_i) \ar[r] \ar[dd] & 
\prod_i L(N_m - V_i) \ar[dd] \\
& \\
L(N_m - N_m^*)  \ar[r] & 
L(N_m) 
\\}
\] 
Here, we have total
degree zero for $L(N_m)$, hence degree $-m$ for $L(N_m^*)$.
We thus get a canonical morphism of complexes of Nisnevich sheaves
\[
L(\FN') \longto L(\FN) \; .
\]
The Mayer--Vietoris property
\cite[Prop.~V.3.1.3]{VSF} shows that
it induces an isomorphism in $\DM$. 
By Lemma~\ref{8B}, the projection from $L(\FN')$ to its upper line 
\[
L(\FN') \longto L(N_m^*) [m]
\]
induces an isomorphism in $\DM$ as well. 

(b) results from (a) by dualizing.
\end{Proofof}

As the proof of Theorem~\ref{8main}
shows, the dependence of the isomorphism on the order $\prec$ can be made explicit.
A canonical isomorphism between $\Mgm ((\doX)^m, \iXnmjZ)$ and a twisted form of $\Mgm (N_m^*) [m]$
can be constructed by proceeding as in \cite[(3.1.4)]{D}. We leave the details to the reader.


\bigskip

%
%

\end{document}